\documentclass[12pt]{amsart}

\usepackage{latexsym}
\usepackage[all]{xy}
\usepackage{amsmath, amsthm}
\usepackage{amssymb}
\usepackage{amsfonts}
\usepackage{color}
\usepackage{mathtools}
\usepackage[T1]{fontenc}

\usepackage{hyperref}
\usepackage{xcolor}
\hypersetup{colorlinks=true,urlcolor=blue,linkcolor=blue,citecolor=blue}

\usepackage{cite}

\newcommand {\Tra} {\widetilde{\mathsf{Traj}}}

\newcommand {\Eu} {\widetilde{Eu}}
\newcommand {\Sing} {\mathsf{Sing}}
\newcommand {\Perf} {\mathsf{Perf}}
\newcommand {\QCoh} {\mathsf{QCoh}}
\newcommand {\Coh} {\mathsf{Coh}}

\newcommand {\Map} {\mathbf{Map}}
\newcommand {\Parf} {\mathsf{Perf}}
\newcommand {\rh} {\mathbb{R}\underline{Hom}}
\newcommand {\rch} {\mathbb{R}\underline{\mathcal{H}om}}
\newcommand {\OO} {\mathcal{O}}

\newcommand {\F} {\mathcal{F}}

\newcommand {\A} {\mathcal{A}}
\newcommand {\M} {\mathcal{M}}
\newcommand {\T} {\mathbb{T}}

\newcommand{\QQ}{\mathbb{Q}}
\newcommand{\VV}{\mathbb{V}}

\newcommand{\ZZ}{\mathbb{Z}}
\newcommand{\LL}{\mathbb{L}}

\newcommand {\Spec} {\mathbf{Spec}}

\newcommand  {\dSt}   {\mathbf{dSt}}

\newcommand{\s}{\infty}
\newcommand{\HH}{\mathbb{H}}

\newcommand{\E}{\mathbb{E}}

\newcommand{\Mgn}{\bar{\mathcal{M}}}
\newcommand{\RMgn}{\mathbb{R}\bar{\mathcal{M}}}
\newcommand{\Cgn}{\bar{\mathcal{C}}}

\newcommand{\D}{\mathcal{D}}
\newcommand{\B}{\mathcal{B}}

\newcommand{\Gm}{\mathbb{G}_m}
\newcommand{\Ql}{\mathbb{Q}_{\ell}}

\theoremstyle{plain}
\newtheorem{thm}{Theorem}[subsection]
\newtheorem{df}[thm]{Definition}
\newtheorem{prop}[thm]{Propositon}
\newtheorem{rmk}[thm]{Remark}

\newtheorem{lem}[thm]{Lemma}

\begin{document}

\title{Foliations and stable maps}
\author{Bertrand To\"en and Gabriele Vezzosi}


\maketitle

\begin{abstract} 
This paper is part of an ongoing series of works on the study of foliations on algebraic varieties 
via derived algebraic geometry. We focus here on the specific case of globally defined vector fields and the global behaviour of their algebraic integral curves. For a smooth and proper variety $X$ with a global vector field $\nu$, we consider the induced vector field $\nu_{g,n}$ on the derived stack of stable maps, of genus $g$ with $n$ marked points, to $X$. When $(g,n)$ is either $(0,2)$ or $(1,0)$, the derived stack of zeros of $\nu_{g,n}$ defines a proper \emph{moduli of algebraic trajectories} of $\nu$. When $(g,n)=(0,2)$ algebraic trajectories behave
very much like rational algebraic paths from one zero of $\nu$ to another, and in particular they can be composed. 
This composition is represented by the usual gluing maps in Gromov-Witten theory, and we use it give three categorical constructions, of different categorical levels, related, in a certain sense, by decategorification. In order to do this, in particular, we have to deal with virtual fundamental classes of non-quasi-smooth derived stacks. When $(g,n)=(1,0)$, zeros of 
$\nu_{1,0}$ might be thought as algebraic analogues of periodic orbits of vector fields on smooth real manifolds. In particular, we propose a Zeta function counting the zeros of $\nu_{1,0}$,  that we like to think of as 
an algebraic version of Ruelle's dynamical Zeta function. We conclude the paper with a brief indication on how to extend these results to the case of general one dimensional foliation $F$, by considering the derived stack of $F$-equivariant stable maps.

\end{abstract}

\tableofcontents

\section*{Introduction}

This paper is part of an ongoing series of works on the study of foliations on algebraic varieties 
via derived techniques (\cite{tvrh, tvindex, tvbook}). We focus here on the specific case of globally defined vector fields (and more generally one dimensional foliations, see our section \S 4) and the global behaviour of their algebraic integral curves. 
The results presented in this paper are preliminary, and we will limit ourselves to survey 
the main constructions and results, more details as well as more general results will appear elsewhere. \\

For a smooth and proper variety $X$ and a global vector field $\nu \in H^0(X,\T_X)$, 
we consider curves in $X$ which are globally invariant by $\nu$, with a particular emphasis on their
moduli and on the various numerical and homological invariants they give rise to.\\ To start with, we will show how $\nu$ 
defines an induced vector field $\nu_{g,n}$ on $\RMgn_{g,n}(X)$, the derived moduli stack of 
stable maps of genus $g$ with $n$ marked points in $X$ (see Section \ref{recoll}). 
An important first observation is that when $(g,n)$ is either $(0,2)$ or $(1,0)$, the zeros of the vector field $\nu_{g,n}$
on $\RMgn_{g,n}(X)$ are precisely given by stable maps that are globally invariant by $\nu$, 
or, equivalently, by \emph{algebraic trajectories} of $\nu$ in $X$. The derived stack 
of zeros of $\nu_{g,n}$ therefore defines a proper moduli space of algebraic trajectories of $\nu$, 
which is the central object of this work. 

When $(g,n)=(0,2)$, algebraic trajectories appear as \emph{broken rational curves} in $X$ along $\nu$
whose marked points are sent to zeros of the original vector field $\nu$ on $X$. Therefore, they behave
very much like algebraic paths from one zero of $\nu$ to another, and in particular they can be composed. 
This composition is represented by gluing maps on the level of derived moduli stacks, induced by the usual gluing 
construction of stable maps in Gromov-Witten theory. Out of this, we extract three constructions, of
different categorical levels, that can be thought as the starting point of an \emph{enumerative geometry of algebraic
trajectories}. We begin by considering derived categories attached to the derived moduli $Z$ of zeros
of $\nu_{0,2}$, namely the derived category of coherent complexes on $Z$ which are relatively perfect along 
the evaluation map $Z \to Z_0\times Z_0$ where $Z_0$ is the derived zero locus of $\nu$ (see Definition \ref{d7}). As base change always holds in derived algebraic geometry (see \cite[3.1]{toenems}), 
the gluing maps induce a composition on these derived categories, and we get this way a linear $2$-category
$\Tra(X,\nu)$ whose set of objects is the set of zeros of $\nu$ in $X$, and whose categories of maps are derived
categories of the derived stack of zeros of $\nu_{0,2}$ (see Definition \ref{d7}). This structure can be viewed  as a
generalization of the categorical Gromov-Witten invariants of \cite{maro} suitably localized around the zeros of the vector field $\nu$. It is also 
a vector field analogue of the usual localization techniques in Gromov-Witten theory for $\Gm$-actions (see for instance \cite{Graber}). 

The $2$-category $\Tra(X,\nu)$ can then be decategorified twice. As a first decategorification,
we construct (Definition \ref{d8}) a $1$-category, whose set of objects is again the zeros of $\nu$, and the
space of morphisms are given by Borel-Moore homology of the derived zero loci of $\nu_{0,2}$ (again relative
to the evaluation map). A second decategorification is obtained by 
considering virtual fundamental classes of the derived zero loci of $\nu_{0,2}$. This is now
a bivector in the homology of the zero locus of $\nu$: a vector field
version of the usual Gromov-Witten invariants in cohomology. On a technical level, the derived moduli stacks involved in this
construction are not quasi-smooth in general, and thus their virtual fundamental classes do not 
exist in the classical sense of the word. We show however that one can always define them as \emph{rational functions} in an indeterminate $t$, 
and that the bivector mentioned above has only $\QQ(t)$ coefficients. More precisely, if $Z_0$ is the zero locus of $\nu$ in $X$, our construction 
provides an element in $H_*(Z_0)^{\otimes 2}\otimes \QQ(t)$.
When further assumptions on $X$ are made (
for instance when $X$ is convex), the derived moduli spaces are indeed quasi-smooth, and the above rational functions are
polynomials which can then be evaluated at $t=-1$. Under these conditions, we obtain a linear 
endomorphism of the vector space spanned by zeros of $\nu$ (assuming these zeros are simple), 
which we like to think of an algebraic version of the \emph{Morse differential}. Finally, we prove (Proposition \ref{p5}) that 
this endomorphism, even if it does not square to zero, satisfies a remarkable formula close to 
the splitting/associativity formula in usual Gromov-Witten theory (see e.g. \cite[Def. 7.1-(3)]{Manin}).

In Section \ref{sectionzeta} we share some thoughts on the case $(g,n)=(1,0)$, where zeros of 
$\nu_{1,0}$ can be thought as algebraic analogues of periodic orbits of vector fields on smooth real manifolds. In particular, we
use our numerical invariants with coefficients in $\QQ(t)$ in order to 
propose a formula for a \emph{Zeta function} counting zeros of $\nu_{1,0}$, which we like to think of as 
an algebraic version of the dynamical Zeta function of a flow (see e.g.
\cite[\S 1.4-(B)]{Ruelle}). We hope to come back to this aspect in future work. 
We conclude the paper with a short section (Section \ref{sectionhigher}) explaining how to extend our results to the case of general one dimensional foliations
instead of just global vector fields, which will certainly provide a wider range of applicability of our results. \\ We recall in the Appendix the basic facts we need about the $K$-theory and $G$-theory of quasi-smooth derived stacks.\\

\textbf{Conventions.} We will work over base field $k$ of characteristic $0$. 
However, all our results and construction remains correct over an arbitrary 
excellent base ring $k$ (but some algebraic operations such as 
symmetric dg-algebras have to be replaced with their simplicial 
versions). 

We will work in the context of derived algebraic geometry of
\cite{hagII,toenems}\\

\section{Vector fields and stable maps}

In this first section we show how a global vector field on an algebraic 
variety $X$ induces
a vector field on $\RMgn_{g,n}(X)$, the derived moduli stack of 
stable maps to $X$. When $(g,n)=(0,2)$ or $(g,n)=(1,0)$, 
we explain how the zeros of the induced vector field correspond to 
\emph{algebraic trajectories} of the original vector field, namely 
stable maps which are globally invariant.

\subsection{Recollections on the derived moduli of stable maps}\label{recoll}

The reader can find details concerning the derived moduli stack of stable
maps in \cite{MR3341464}. We remind below some notations and constructions.

We let $\Mgn_{g,n}^{pr}$ be the Artin stack of prestable curves of genus $g$ with $n$ marked points. 
We denote by $\pi : \Cgn_{g,n} \to \Mgn_{g,n}^{pr}$ be the universal prestable curve. The morphism
$\pi$ induces an adjunction on the $\s$-categories of relative derived stacks
$$\pi^* : \dSt/\Mgn_{g,n}^{pr} \leftrightarrows \dSt/\Cgn_{g,n} : \pi_*,$$
where $\pi^*$ is the usual base change $\s$-functor. The right adjoint $\pi_*$ 
is also known as the Weil restriction along $\pi$. 

Because $\pi$ is a proper and flat morphism, it is known that $\pi_*$ preserves derived Artin stacks locally
of finite presentation. 
In particular, if $X$ is any smooth and proper scheme over $k$, we can consider
$$\RMgn_{g,n}^{pr}(X):=\pi_*(X \times \Cgn_{g,n}),$$
which is a derived Artin stack endowed with a natural projection to $\RMgn_{g,n}^{pr}$.

\begin{df}\label{d1}
The \emph{derived stack of pre-stable maps on $X$} (of genus $g$ with $n$ marked points)
is the derived Artin stack $\RMgn_{g,n}^{pr}(X)$ defined above.
\end{df}

The projection $p : \RMgn_{g,n}^{pr}(X) \to \Mgn_{g,n}^{pr}$ is a morphism 
of derived Artin stacks, whose fibers are the derived mapping schemes of pre-stable curves
to $X$. In more precise terms, for any derived scheme $S$ and morphism $S \to \Mgn_{g,n}^{pr}$, 
the pull-back $S\times_{\Mgn_{g,n}^{pr}}\RMgn_{g,n}^{pr}(X)$ is canonically identified with
$\Map_{S}(\Cgn_{g,n}\times_{\Mgn_{g,n}^{pr}}S,X\times S)$, the derived mapping stack over $S$
from the curve $\Cgn_{g,n}\times_{\Mgn_{g,n}^{pr}}S$ to $X\times S$. We see in particular that 
the morphism  $p : \RMgn_{g,n}^{pr}(X) \to \Mgn_{g,n}^{pr}$ is representable and its fibers are
derived schemes locally of finite type.

As usual, the derived stack of stable maps is defined as the open substack 
$\RMgn_{g,n}(X) \subset \RMgn_{g,n}^{pr}(X)$ whose geometric points consists of 
stable maps to $X$. The derived stack $\RMgn_{g,n}(X)$ is Deligne-Mumford, and 
each component $\RMgn_{g,n}(X,\beta)$, of stable maps $f : C \to X$ 
with $f_*([C])=\beta$, is proper over $\Spec\, k$ (see \cite[Thm. 3.14]{bema}). 
Here $\beta$ is an element in $H_2(\overline{X},\Ql(-1))$, 
is a class in $\Ql$-adic homology of $\overline{X}=X \otimes_k \overline{k}$.

We then define the open (and closed) substacks $\RMgn_{g,n}(X,\beta) \subset 
\RMgn_{g,n}(X)$ by considering geometric points $f : \Spec\, K \to \RMgn_{g,n}(X)$
corresponding to stable maps $C \to X\otimes_k \overline{k}$, 
with $C$ a prestable curve over 
$K$, such that $f_*([C])=\beta \in H_2(X\otimes_k K,\Ql(-1))$, where
$[C] \in H_2(C,\Ql(-1))$ is the fundamental class of $C$.

\begin{df}\label{d2}
The \emph{derived stack of stable maps on $X$} (of genus $g$ with $n$ marked points)
is the derived Deligne-Mumford stack $\RMgn_{g,n}(X)$ defined above. The open and 
closed
substack of stable maps of class 
$\beta \in H_2(X,\Ql(1))$ is $\RMgn_{g,n}(X,\beta)$ defined above.
\end{df}

An important fact that we will use in a crucial way is the
explicit computation of the tangent complex of $\RMgn_{g,n}^{pr}(X)$. For this, we 
consider the evaluation morphism
$$ev : \RMgn_{g,n}^{pr}(X) \times_{\Mgn_{g,n}^{pr}}\Cgn_{g,n} \longrightarrow X,$$
as well as the projection to the first factor
$$q : \RMgn_{g,n}^{pr}(X) \times_{\Mgn_{g,n}^{pr}}\Cgn_{g,n} \longrightarrow \RMgn_{g,n}^{pr}(X).$$
The tangent complex of $p : \RMgn_{g,n}^{pr}(X) \to \Mgn_{g,n}^{pr}$ is then given by (see \cite[\S 3.4]{toenems})
$$\T_{\RMgn_{g,n}^{pr}(X)/\Mgn_{g,n}^{pr}} \simeq q_*(ev^*(\T_X)),$$
where $T_X$ is the tangent vector bundle of $X$ relative to $k$. As $q$ is a flat and proper curve, 
we see immediately that $\T_{\RMgn_{g,n}^{pr}(X)/\Mgn_{g,n}^{pr}}$ is a perfect complex whose amplitude is contained
in $[0,1]$. To obtain the full tangent complex of $\RMgn_{g,n}^{pr}(X)$, relative to $k$, we use the exact triangle
$$\xymatrix{\T_{\RMgn_{g,n}^{pr}(X)/\Mgn_{g,n}^{pr}} \ar[r] & \T_{\RMgn_{g,n}^{pr}(X)} \ar[r] & 
p^*(\T_{\Mgn_{g,n}^{pr}}).}$$
As $\Mgn_{g,n}^{pr}$ is a smooth Artin stack over $k$, $p^*(\T_{\Mgn_{g,n}^{pr}})$ is itself perfect 
and of amplitude in $[-1,0]$. As a result, the tangent complex $\T_{\RMgn_{g,n}^{pr}(X)}$ is
perfect and of amplitude in $[-1,1]$. The derived Artin stack 
$\RMgn_{g,n}^{pr}(X)$ is thus quasi-smooth over $k$ (see Appendix), and therefore
$\RMgn_{g,n}(X)$ is itself a quasi-smooth and proper derived Deligne-Mumford stack over $k$
(its tangent complex is a complex and of Tor-amplitude $[0,1]$).

\subsection{Induced vector fields}

We keep working with $X$ a smooth and proper scheme over $k$. As before $\T_X$ denotes the tangent bundle
of $X$ relative to $k$.

\begin{prop}\label{p1}
Any vector field $\nu \in H^0(X,\T_X)$ determines a natural vector field
$$\nu_{g,n} \in H^0(\RMgn_{g,n}(X),\T_{\RMgn_{g,n}(X)}).$$
\end{prop}

\textit{Proof.} Start with $\nu : \OO_X \to \T_X$ on $X$. By pull-back along the evaluation map
$ev : \RMgn_{g,n}^{pr}(X) \times_{\Mgn_{g,n}^{pr}}\Cgn_{g,n} \longrightarrow X$, it provides
a new section $ev^*(\nu) : \OO_{\RMgn_{g,n}^{pr}(X)} \to ev^*(\T_X)$. Its push-foward by 
$q : \RMgn_{g,n}^{pr}(X) \to \Mgn_{g,n}^{pr}$ produces a morphism
$$\OO_{\RMgn_{g,n}^{pr}} \longrightarrow q_*q^*(\OO) \longrightarrow q_*ev^*(\T_X).$$
We now use that $q_*ev^*(\T_X) \simeq \T_{\RMgn_{g,n}^{pr}(X)/\Mgn_{g,n}^{pr}}$, as well as
the natural morphism $\T_{\RMgn_{g,n}^{pr}(X)/\Mgn_{g,n}^{pr}} \to \T_{\RMgn_{g,n}^{pr}(X)}$
in order to get a canonical element $\nu_{g,n} \in H^0(\RMgn_{g,n}(X),\T_{\RMgn_{g,n}(X)})$. \hfill $\Box$ \\

\begin{rmk}\emph{One can understand the above proposition in various ways. For instance, it can 
be seen as the fact that an infinitesimal automorphism of $X$ induces
an infinitesimal automorphism of $\RMgn_{g,n}^{pr}(X)$, simply by using the functoriality 
of the construction $X \mapsto \RMgn_{g,n}^{pr}(X)$. In another direction, and more generally, a global vector
field defines a derived foliation of dimension $1$ in the sense of \cite{tvrh}, and
the above proposition is a very particular case of a more general result stating the existence
of direct images of derived foliations (see \cite{tvbook}).}
\end{rmk}

It will be important to remind (see the proof of Proposition \ref{p1}) that the vector field $\nu_{g,n}$ is in fact the image
of a natural relative vector field $\nu'_{g,n} \in H^0(\RMgn_{g,n}^{pr}(X),\T_{\RMgn_{g,n}^{pr}(X)/\Mgn_{g,n}^{pr}})$
by the canonical map $\T_{\RMgn_{g,n}^{pr}(X)/\Mgn_{g,n}^{pr}} \to \T_{\RMgn_{g,n}^{pr}(X)}$, and thus is vertical relative to this morphism.
Our main object of interest is $\nu_{g,n}$, but it will be sometimes important to use $\nu'_{g,n}$ in some
of the constructions. We thus gather these two objects $\nu_{g,n}$ and $\nu_{g,n}'$ is the definition below.

\begin{df}\label{d3}
For $\nu \in H^0(X,\T_X)$ a global vector field on $X$ (relative to $k$). The \emph{induced vector field}
(resp. \emph{relative induced vector field}) \emph{of type $(g,n)$} is $\nu_{g,n}$ (resp. $\nu'_{g,n}$) constructed above.
\end{df}

Pointwise, the vector field $\nu_{g,n}$ can be described as follows. 
Let $f : C \to X$ by a stable map to $X$ (say define over $k$), 
and $\Sigma \subset C$ be the divisor defining the marked points on $C$.
We consider $\T_C(-\Sigma)$, the tangent complex of $C$ twisted by $\OO(-\Sigma)$, 
and $\T_C(-\Sigma) \to f^*(\T_X)$ the induced by the derivative of $f$. The cone
of this morphism is a perfect complex $E$ on $C$ of amplitude 
contained in $[-1,0]$ on $C$.
The tangent complex of $\RMgn_{g,n}(X)$ at the point $C$ is given by the complex
of $k$-vector spaces
$H^*(C,E)$.

We have a canonical morphism $f^*(\T_X) \to E$, as well as
a global section $\OO_C \to f^*(\T_X)$ given by the restriction
of $\nu$ on $C$. The vector field $\nu_{g,n}$
evaluated at the point $C$, is then given by the composition
$$k \to H^*(C,\OO_C) \to H^*(C,f^*(\T_X)) \to H^*(C,E).$$

Let $X$ be a smooth and proper scheme over $k$ and $\nu \in H^*(X,\T_X)$. 
The zeros of the induced vector $\nu_{g,n}$ on $\RMgn_{g,n}(X)$ defines a closed derived substack
$Z(\nu_{g,n}) \subset \RMgn_{g,n}(X)$. More precisely, 
we can form the sheaf of commutative simplicial rings $Sym_{\OO_{\RMgn_{g,n}(X)}}(\LL_{\RMgn_{g,n}(X)})$, 
whose spectrum
(relative to $\RMgn_{g,n}(X)$)
is the total tangent stack 
$$\T\RMgn_{g,n}(X)=\Spec_{\RMgn_{g,n}(X)}\, Sym_{\OO_{\RMgn_{g,n}(X)}}(\LL_{\RMgn_{g,n}(X)}).$$
The augmentation $Sym_{\OO_{\RMgn_{g,n}(X)}}(\LL_{\RMgn_{g,n}(X)}) \to \OO_{\RMgn_{g,n}(X)}$ defines a natural projection
$\T\RMgn_{g,n}(X) \to \RMgn_{g,n}(X)$.
The vector field $\nu_{g,n}$ provides a natural section 
$s : \RMgn_{g,n}(X) \to \T\RMgn_{g,n}(X)$. Similarly, the zero vector field provides
a section $0 : \RMgn_{g,n}(X) \to \T\RMgn_{g,n}(X)$. The derived stack of zeros
of $\nu_{g,n}$ is naturally defined as the derived intersections of these two sections.

\begin{df}\label{d4}
With the above notations, the \emph{derived stack of zeros of $\nu_{g,n}$}
is defined by 
$$Z(\nu_{g,n}):= \RMgn_{g,n}(X) \times_{s,\T\RMgn_{g,n}(X),0} \RMgn_{g,n}(X).$$
\end{df}

The derived stack $Z(\nu_{g,n})$ can be seen as a closed substack of $\Mgn_{g,n}(X)$ my means of the
second projection 
$$j : Z(\nu_{g,n}):= \RMgn_{g,n}(X) \times_{s,\T\RMgn_{g,n}(X),0} \RMgn_{g,n}(X) \longrightarrow 
\RMgn_{g,n}(X).$$

The derived substack $Z(\nu_{g,n})$ can also be described in terms of a relative spectrum 
of an explicit sheaf of commutative simplicial rings, as follows. The section $\nu_{g,n}$ defines
a morphism of perfect complexes on $\RMgn_{g,n}(X)$
$$\LL_{\RMgn_{g,n}(X)} \longrightarrow \OO_{\RMgn_{g,n}(X)}.$$
By the universal properties of symmetric algebras, this morphism defines a morphism
of sheaves of commutative simplicial rings
$$Sym_{\OO_{\RMgn_{g,n}(X)}}(\LL_{\RMgn_{g,n}(X)}) \longrightarrow \OO_{\RMgn_{g,n}(X)}.$$
In the same manner we have the canonical augmentation morphism
$$Sym_{\OO_{\RMgn_{g,n}(X)}}(\LL_{\RMgn_{g,n}(X)}) \longrightarrow \OO_{\RMgn_{g,n}(X)}.$$
Then, the derived substack $Z(\nu_{g,n})$ can also be described as
$$Z(\nu_{g,n}) \simeq \Spec\, \left( \OO_{\RMgn_{g,n}(X)} \otimes_{Sym_{\OO_{\RMgn_{g,n}(X)}}(\LL_{\RMgn_{g,n}(X)})} 
\OO_{\RMgn_{g,n}(X)} \right).$$

A first important comment is the fact that, even though $\RMgn_{g,n}(X)$ is quasi-smooth, 
the derived stack $Z(\nu_{g,n})$ is no more quasi-smooth in general. Indeed, 
its tangent complex is in general concentrated in amplitude $[0,2]$. As an example, 
when $\nu=0$ is the zero vector field on $X$, then $Z(\nu_{g,n}) \simeq \T\Mgn_{g,n}(X)[-1]$
is the shifted total tangent space, whose tangent complex taken at a geometric point $x$
is $\T_{\Mgn_{g,n}(X),x} \oplus \T_{\Mgn_{g,n}(X),x}[-1]$. As a consequence, 
the derived stack $Z(\nu_{g,n})$ does not carry a virtual class (see Appendix) as its structure sheaf 
is in general not bounded. We will see however in the next sections that it is possible
to overcome this problem and to associate to $Z(\nu_{g,n})$ a homology class that 
deserves the name of \emph{virtual class} and that can be used in order to make several
interesting constructions.

As a second comment, a global point in $Z(\nu_{g,n})$ morally corresponds to 
a stable map $f : C \to X$ which is globally preserved by the vector field $\nu$. 
In particular, for such point, $\nu$ defines a global section 
in $H^0(C,\mathbb{T}_C(-\Sigma))$, thanks to the short exact sequence
$$0\to H^0(C,\mathbb{T}_C(-\Sigma)) \to H^0(C,f^*(\mathbb{T}_X)) \to 
H^0(C,E).$$

As a result, when $C$ is a stable curve,
then $H^0(C,\mathbb{T}_C(-\Sigma))\simeq 0$ and thus the map $f$ globally
factors through the zero locus of $\nu$ in $X$. This shows that 
the possibly interesting cases only appear for stable maps from unstable curves
and thus for curves with components of genus $0$ or $1$. We will therefore focus on 
genus $0$ and $1$ in the sequel.

\subsection{Genus $0$ and $1$}

Let $f : \Spec\, L \to \RMgn_{g,n}(X)$ be a morphism with $L$ a separably closed field. By replacing $k$ by $L$ and
$X$ by $X\otimes_{k} L$, we can furthermore assume that $k=L$. The morphism $f$ is then 
described as a stable map
$$f : C \longrightarrow X,$$
with $C$ a prestable curve of type $(0,2)$ or $(1,0)$. \\

\noindent \textsf{Case $(g,n)=(0,2)$.} Let us first assume that $(g,n)=(0,2)$ and that $C$ is smooth to make things simplier, so it consists
of a $\mathbb{P}^1$ marked at $0$ and $\s$. We can describe the vector field $\nu_{g,n}$
at the point $f$ as follows. The tangent complex of $\RMgn_{g,n}(X)$ at $f$ is given by
(see \cite{MR3341464})
$$\T_{\RMgn_{g,n}(X),f}\simeq \HH(C,\T_C(-0-\s) \to f^*(\T_X)).$$
The vector field $\nu : \OO_X \to \T_X$ provides a morphism $\HH(C,\OO_C) \to \HH(C,f^*(\T_X)) \to 
\HH(C,\T_C(-0-\s) \to f^*(\T_X))$. Since $C=\mathbb{P}^1$, we have 
$\HH(C,\OO_C) \simeq k$, and the corresponding morphism of complexes of $k$-vector spaces
$$k \longrightarrow \HH(C,\T_C(-0-\s) \to f^*(\T_X))$$
is the vector field $\nu_{0,2}$ evaluated at the point $f$. We also note that 
$k \to \HH(C,f^*(\T_X))$ is $\nu'_{2,0}$ evaluated at $f$, where $\nu'_{0,2}$ is the
relative induced vector field of Definition \ref{d3}.

We now consider the commutative diagram
$$\xymatrix{\HH(C,f^*(\T_X)) \ar[rr] & & \HH(C,\T_C(-0-\s) \to f^*(\T_X)) \\
 & \ar[ul]^-{\nu'_{0,2}}k \ar[ru]_-{\nu_{0,2}} & 
}$$
A homotopy to zero of $\nu_{0,2}$ then consists of a factorization of $\nu'_{0,2}$ via the
canonical morphism induced by the tangent map of $f$.
$$Tf : \HH(C,\T_C(-0-\s)) \longrightarrow \HH(C,f^*(\T_X)).$$
Again as $C=\mathbb{P}^1$, we have $\HH(C,\T_C(-0-\s)) \simeq k$, and a canonical generator
is given by the vector field $z.\partial_z$ on $\mathbb{P}^1$ which vanishes at $0$ and $\s$. 
We thus arrive at the conclusion that the geometric point $f$ lies in the substack $Z(\nu_{0,2})$
if and only if there exists $\lambda \in k$ such that the image of $\lambda.z\partial_z$ by $Tf$
is $\nu$ restricted to $C$. In other words, $f$ lies in $Z(\nu_{0,2})$ if and only if
$f$ exhibits $C$ as an algebraic trajectory of $\nu$ in $X$.

We have seen that a stable map $f : (\mathbb{P}^1,0,\s) \longrightarrow X$, which is a zero
of the induced vector field $\nu_{0,2}$ must be an algebraic trajectory of $\nu$ on $X$. 
In particular, the images of $0$ and $\s$ in $X$ must be a zero of $\nu$ (and moreover these
zeros must be non-degenerate along $C$ itself). In general, when $C$ possesses several
irreducible components, and thus is a concatenation of $\mathbb{P}^1$ glued along $0$ and $\s$, 
$f : C \to X$ lies in $Z(\nu_{0,2})$ if and only if each component $C_i$ of $C$ is an algebraic
trajectory of $\nu$ in $X$. Therefore, $f$ becomes a trajectories of broken $\mathbb{P}^1$ in $X$ for
the vector field $\nu$, each irreducible components of $C$ sitting between two zeros of $\nu$.\\

\noindent \textsf{Case $(g,n)=(1,0)$.} Let us now assume that $(g,n)=(1,0)$. The same reasoning shows that the geometric points of $Z(\nu_{1,0})$
correspond to stable maps $f : E \to X$, with $E$ a curve of genus $1$, such that there is
a global vector field $\alpha \in H^0(E,\T_E) \simeq k$ whose image by $Tf$ is the restriction 
of $\nu$ to $E$. These are again the trajectories of $\nu$ in $X$.\\

As a conclusion, we have defined for any vector field $\nu$ on $X$ an induced vector
field $\nu_{g,n}$ on the derived stack of stable maps $\RMgn_{g,n}(X)$. When $(g,n)$ is either
equal to $(0,2)$, or to $(1,0)$, the zeros of $\nu_{g,n}$ on $\RMgn_{g,n}(X)$ precisely are 
the algebraic trajectories of $\nu$ on $X$. Therefore, the derived stack of zeros
$Z(\nu_{g,n})$ is a Deligne-Mumford derived stack over $k$, which can be written 
as a disjoint union of proper Deligne-Mumford derived stack
$Z(\nu_{g,n}) \simeq \coprod_{\beta} Z(\nu_{g,n})(\beta)$, 
and is the starting point to study enumerative properties of algebraic trajectories.
However, 
the derived stack $Z(\nu_{g,n})$ is not quasi-smooth in general, and 
extracting numerical or cohomological invariants out of it requires more work. We will see 
in the next section that $Z(\nu_{g,n})$ still possesses a nice cohomology theory and
an associated virtual class, which can be used in order to 
consider enumerative questions concerning algebraic trajectories of the original vector field $\nu$. Also, when $X$ is convex $Z(\nu_{g,n})$ is automatically
quasi-smooth and things simplify a lot (see our section \S 3.2).

\section{Localized Euler class in the derived setting}

In this section we study the zero locus of a section $s$ of 
a perfect complex $E$, of amplitude contained in $[0,1]$, over
a quasi-smooth derived Deligne-Mumford stack $X$. We will see how we can 
associate several cohomological, numerical and categorical invariants
to such data. In the next section, it will be applied to the special
case where $E= \T_{\RMgn_{g,n}(X)}$ and $s=\nu_{g,n}$ is the induced vector field of Definition \ref{d3}. 

We fix once for all a quasi-smooth derived Deligne-Mumford stack $X$ over $k$, 
$E \in \Parf(X)$ of amplitude contained in $[0,1]$ and 
$s : \OO_X \to E$ a morphism in $\Parf(X)$. We denote by $Z(s) \subset X$ the derived
substack of zeros of $s$, which is defined to be 
$\Spec\, (\OO_X \otimes_{Sym_{\OO_X}(E^\vee)}\OO_X)$, where the first map
$Sym_{\OO_X}(E^\vee) \to \OO_X$ is induced by the section $s$, and the second
map $Sym_{\OO_X}(E^\vee) \to \OO_X$ is the natural augmentation (i.e. the previous map, with $s=0$). 

When the perfect complex $E$ is of amplitude $0$ (i.e. is a vector bundle), 
 $Z(s)$ is quasi-smooth and
thus gives rise to a virtual fundamental class $[Z(s)]^{vir} \in G_0(Z(s))$, in the Grothendieck
group of coherent sheaves (or in some Borel-Moore homology theory, see Appendix). However, when $E$
is no more a vector bundle, then $Z(s)$ is not quasi-smooth in general, as 
it is the derived intersection of quasi-smooth derived stacks inside an ambient quasi-smooth
derived stack. The virtual class $[Z(s)]^{vir}$ thus no more makes sense in this generality. However,
we will see that the Euler class do exists as a $\Gm$-equivariant class (with support on $Z(s)$)
on $\VV(E[1])=\Spec\, (Sym_{\OO_X}(E[1]))$, the linear stack associated to the shifted complex $E[1]$. 
The purpose of this section is to explain the construction of this class, as well as of the
numerical and cohomological invariants that can be extracted from it. We will also provide a categorical interpretation
in terms of certain categories of singularities associated to functions of degree $-1$ (Section \ref{seccatsing}).

\subsection{The refined Euler class}

We consider the linear derived 
stack $\VV(E^\vee)=\Spec\, Sym_{\OO_X}(E^\vee)$, associated to the dual complex $E^\vee$ (see \cite{monier,toenems} for more on linear stacks). 
It comes equipped with a natural projection $\pi : \VV(E^\vee) \to X$ making it into
an affine derived scheme relative to $X$.
The derived stack $\VV(E^\vee)$ is thus itself 
a derived Deligne-Mumford stack which is affine over $X$. Moreover, 
as $E^\vee$ has its amplitude contained in $[-1,0]$, $Z$ is also quasi-smooth and so is the natural
projection $\pi : \VV(E^\vee) \to X$.

The projection $\pi$ possesses two sections $e : X \to \VV(E^\vee)$ and $s : X \to \VV(E^\vee)$. The section $e$ is the
zero section and corresponds to the natural augmentation $Sym_{\OO_X}(E^\vee) \to \OO_X$. The section 
$s$ is induced by $s : \OO_X \to E$ and by the induced morphism of commutative $\OO_X$-algebras
$Sym_{\OO_X}(E^\vee) \to \OO_X$. We have an $\s$-functor
$$\underline{Hom}_{\OO_X}(s^*e_*(\OO_X),-) : \QCoh(X) \longrightarrow \QCoh(X).$$
We let $\A:=\underline{Hom}_{\OO_{\VV(E^\vee)}}(e_*(\OO_X),e_*(\OO_X))$, considered 
as a quasi-coherent associative $\OO_X$-algebra. By construction, the above $\s$-functor 
factors as
$$\xymatrix{
\QCoh(X) \ar[r] & \A-\mathsf{Mod}^{qcoh} \ar[r] & \QCoh(X),
}$$
where $\A-\mathsf{Mod}^{qcoh}$ is the $\s$-category of 
quasi-coherent $\A$-modules, and the second $\s$-functor is the
forgetful $\s$-functor. The first $\s$-functor will be denoted by
$Kos_s : \QCoh(X) \to \A-\mathsf{Mod}^{qcoh}$.

By construction, the underlying quasi-coherent complex of $\A$
is naturally equivalent to $\oplus_{n\geq 0}(Sym^n_{\OO_X}(E^\vee[1]))^{\vee}$. Therefore, there exists
a canonical inlcusion morphism $E[-1] \to \A$ corresponding to the the $n=1$ factor. This in turn induces
a morphism on cohomology sheaves
$H^1(E)[-2] \to H^*(\A)$. The graded algebra $H^*(\A)$ is moreover graded commutative, as 
$\A$ underlies an $E_\s$-structure, and therefore the above morphism induces a natural morphism of (underived) graded 
$H^0(\OO_X)$-algebras
$$\B:=Sym_{H^0(\OO_X)}(H^1(E)) \longrightarrow H^*(\A),$$
where $H^1(E)$ sits in weight $2$ in the graded algebra $\B$.
We warn the reader here that the $Sym_{H^0(\OO_X)}$ must be understood in the non-derived setting.
In particular, any quasi-coherent $\A$-module $M \in \A-\mathsf{Mod}^{qcoh}$ gives rises
to a sheaf of graded $\B$-modules by considering $H^*(M)$ with its graded $H^*(\A)$-module structure
and induced graded $\B$-module structure from the above morphism. Note that
graded modules over the algebra $\B$ already appears in the construction 
of the singular-support in \cite{arga}[\S 4.1.2].

The following finiteness lemma somehow already appears in 
\cite{arga}[Appendix D]. 
We reproduce a proof here for sake of completeness.

\begin{lem}\label{l1}
For any perfect complex $M \in \Parf(X)$, 
the sheaf of graded $\B$-modules $H^*(Kos_s(M))$ is
locally of finite type.
\end{lem}

\textit{Proof.} The statement is local for the \'etale topology on $X$ and we can thus
assume that $X=\Spec\, A$ is affine. Moreover, as $\OO_X$
generates perfect complexes (because $X$ is assumed to be affine), it is enough to prove the statement for 
$M=\OO_X$. We can present $E$ as the fiber of a morphism $d$ of vector bundles 
$\xymatrix{E \ar[r] & E_0 \ar[r]^{d} & E_1}$. This defines a closed embedding 
of total spaces $j : \VV(E^\vee) \hookrightarrow \VV(E_0^\vee)$ as well as an 
epimorphism of graded algebras $H^0(Sym_A(E_1))\simeq Sym_{H^0(\OO_X)}(H^0(E_1)) \to \B$. 

The derived scheme $\VV(E^\vee)$ sits inside $\VV(E_0^\vee)$ as the fiber at $0$
of the morphism $\VV(E_0^\vee) \to \VV(E_1^\vee)$ induced by $d$. This implies that 
for two bounded coherent complexes $\F$ and $\F'$ on 
$\VV(\E^\vee)$, the quasi-coherent complex $\underline{Hom}(j_*(\F),j_*(\F'))$
comes equipped with a canonical structure of a module over $Sym_{A}(E_1^\vee[1])$. 
This structure is moreover such that we have
$$\underline{Hom}(\F,\F') \simeq \underline{Hom}_{Sym_A(E_1^\vee[1])}(A,\underline{Hom}(j_*(\F),j_*(\F'))).$$
In our situation, we apply this to $\F=e_*(\OO_X)$ and $\F'=s_*(M)$, and we
therefore see that because $M$ is perfect, $j_*(\F)$ and $j_*(\F')$ are again 
perfect over $\VV(E_0^\vee)$. 
Therefore, $\underline{Hom}(j_*(\F),j_*(\F'))$ is perfect $A$-module, and thus it is bounded
with coherent cohomologies. As an $Sym_A(E^\vee[1])$-module, it thus 
lies in the thick triangulated subcategory of generated by $H^0(A)$-module of finite type (considered
as $Sym_A(E^\vee[1])$-modules via the augmentation $Sym_A(E^\vee[1]) \to A \to H^0(A)$). By the formula above,
this implies that $\underline{Hom}(\F,\F')$ is, as a module over 
$\underline{End}_{Sym_A(E_1^\vee[1])}(A)$, in the thick triangulated subcategory generated
by $\underline{End}_{Sym_A(E_1^\vee[1])}(A) \otimes_A P$ with $P$ an $H^0(A)$-module of finite type.
In particular, 
$H^*(\underline{Hom}(\F,\F'))$ is a finitely generated graded module over 
$H^*(\underline{End}_{Sym_A(E_1^\vee[1])}(A))\simeq Sym_{H^0(A)}(H^0(E_1))\otimes_{H^0(A)}H^*(A)$, and thus
over $Sym_{H^0(A)}(H^0(E_1))$ because $H^*(A)$ is a finite $H^0(A)$-module.
This graded module is the image of $H^*(Kos_s(M))$ by the canonical epimorphism 
$Sym_{H^0(\OO_X)}(H^0(E_1)) \to \B$, and thus we conclude that $H^*(Kos_s(M))$ is also finitely generated
as a graded $\B$-module. \hfill $\Box$ \\

We now consider the determinant line bundle $\omega_E = det(E)$ on $X$ and consider 
$Kos_s(\omega_E^{-1}) \in \A-\mathsf{Mod}^{qcoh}$. By the above lemma, we have a graded
$\B$-module $H^*(Kos_s(\omega_E^{-1}))$ of finite type. By construction, 
$\B$ can be canonically identified with $H^0(\OO_{\VV(E[1])})$, the $H^0$-sheaf of rings
of functions on $\VV(E^\vee[1]) \to X$, where the graduation is now induced by the
canonical $\Gm$-action on $\VV(E[1])$ given by weight $2$ dilation along the fibers. 
In order to make this weight $2$-action specific, and to distinguish it 
from the natural weight $1$-action we will consider later, we will use the notation
$\VV^{(2)}(E[1])$.

Thanks to lemma \ref{l1}, the graded
$\B$-module $H^*(Kos_s(\omega_E^{-1}))$
thus defines a class in the Grothendieck
group of graded modules of finite type, which by definition is the $\Gm$-equivariant
Grothendieck group of coherent sheaves on the truncation $\tau_0\VV^{(2)}(E[1])$.
As the pushforward along the inclusion $\tau_0\VV^{(2)}(E[1]) \hookrightarrow \VV^{(2)}(E[1])$ induces an 
isomorphism on Grothendieck group of coherent sheaves (by d\'evissage, see Appendix), 
we get this way a well defined class
$$[H^*(Kos_s(\omega_E^{-1}))] \in G_0^{\Gm}(\VV^{(2)}(E[1])).$$

Finally, we note that $H^i(Kos_s(\omega_E^{-1})) = Ext^i(s^*e_*(\OO_X),\omega_E^{-1})$ is zero
outside of $Z=Z(s) \subset X$ the zero locus of $s$, and thus the class
$[H^*(Kos_s(\omega_E^{-1}))]$ above in fact lies in the $G$-theory with support on $Z$, 
or, again by d\'evissage, in $G_0^{\Gm}(\VV^{(2)}(E[1])_Z)$, where 
$\VV^{(2)}(E[1])_Z = \VV^{(2)}(E[1])\times_X Z$
is the restriction of $\VV^{(2)}(E[1])$ to $Z$.

\begin{df}\label{d5}
The class $[H^*(Kos_s(\omega_E^{-1}))] \in G_0^{\Gm}(\VV^{(2)}(E[1])_Z)$ defined above
is called the \emph{refined Euler class of the section $s$}. It is denoted by 
$$\Eu_E(s) := [H^*(Kos_s(\omega_E^{-1}))] \in G_0^{\Gm}(\VV^{(2)}(E[1])_Z).$$
\end{df}

We conclude this section with a couple of comments. \\
 
\noindent \textbf{The bundle case.} When $E$ is a vector bundle, then 
the graded ring $\B$ is simply $\OO_X$, and $H^*(Kos_s(\omega_E^{-1}))$ is 
a graded $H^0(\OO_X)$-module of finite type. In fact, $Kos_s(\omega_E^{-1})$ is 
a perfect complex on $X$, which by Grothendieck duality is
$s^*e_*(\OO_X)[-r]$, where $r$ is the rank of $E$ (note that $\omega_E[r]$ is the
dualizing complex of the closed embedding $s : X \to \VV(E^\vee)$). Therefore, the
class $[\Eu(s)] \in G_0^{\Gm}(\VV^{(2)}(E[1])_Z)\simeq 
G_0^{\Gm}(Z)\simeq G_0(Z)[t,t^{-1}]$ is nothing else than 
$$\Eu_E(s)=(-1)^r\sum_{i \in \ZZ} [Tor^{\OO_{\VV(E^\vee)}}_i(e_*(\OO_X)\otimes_{\OO_X}s_*(\OO_X))]t^i.$$
The value at $t=-1$ of this Laurent polynomial is, up to a sign, thus the usual 
Euler class of $s$
$$\Eu_E(s)(t=-1) = (-1)^r[s^*e_*(\OO_X)] \in G_0(X)_Z \simeq G_0(Z)$$
where $Z=Z(s)$ is the zero locus of $s$ and $G_0(X)_Z$ the
$G$-theory with support in $Z$. Using base change for the cartesian diagram of derived stacks
$$\xymatrix{Z \ar[r] \ar[d] & X \ar[d]^-{s} \\
X \ar[r]_{e} & \VV(E^\vee),}$$
we see that $\Eu(s)(t=-1)$ is the usual (virtual) fundamental class $[Z]^{vir}:=[\OO_Z] \in G_0(Z)$ 
of the derived scheme $Z$. \\

\noindent \textbf{Interpretation in terms of singular supports.} 
The class $\Eu_E(s)$ of Definition \ref{d5}
is inspired, and very closely related to, the singular support of coherent sheaves of \cite{arga}. 
Indeed, we can apply \cite[Thm. 4.18]{arga} to the bounded coherent complexes $e_*(\OO_X)$ and 
$s_*(\omega_E^{-1})$ on 
the quasi-smooth derived stack $\VV^{(2)}(E^\vee)$, and get 
a graded coherent sheaf on the singular schemes of 
$\VV^{(2)}(E^\vee)$. In our setting, we consider
the relative version over $X$, as both $e_*(\OO_X)$ and $s_*(\omega_E^{-1})$ are relatively 
perfect over $X$, 
and thus this support lies on the relative singular scheme of $\VV(E^\vee) \to X$. This relative 
singular scheme is nothing else than $\VV^{(2)}(E[1])$, and our class $\Eu_E(s)$ is thus an 
incarnation of
the relative singular support of $Hom(e_*(\OO_X),s_*(\omega_E^{-1}))$, or rather its relative 
characteristic cycle (i.e. the singular support but counted with its natural multiplicities). \\

\subsection{The virtual Euler series and Euler numbers}\label{eulseriessec}

In the previous paragraph we have constructed the refined Euler class $\Eu(s)$. We now explain 
how this
class can be somehow pushed down to $Z \subset X$ itself, under the form of a rational function
with coefficients in Borel-Moore homology of $Z$. When $Z$ is proper (e.g. when $X$ is so), 
this can be used in order to define Euler numbers which morally count the number of zeros of the 
section $s$. \\

We start by consider the projection $p : \VV^{(2)}(E[1])_Z \to Z$, pull-back from the projection
$\VV^{(2)}(E[1]) \to X$ along the closed embedding $Z=Z(s) \hookrightarrow X$. The projection 
$p$ is an affine morphism equipped with a natural $\Gm$-action of dilations along the 
fiber (with weight
$2$ in our setting). Algebraically,
this corresponds to the graded simplicial commutative $\OO_Z$-algebra $Sym_{\OO_Z}(E_{|Z}[1])$,
where $E$ 
is given weight $2$. We thus get a direct image
on quasi-coherent complexes
$$p_* : \QCoh([\VV^{(2)}(E[1])_Z/\Gm]) \to \QCoh(Z \times B\Gm).$$
As usual, 
we will identify the $\s$-category $\QCoh(Z \times B\Gm)$ with the $\s$-category of
graded quasi-coherent complexes on $Z$ (see \cite{Moulinos}).
Using that $E$ is perfect, it is easy to see that this push-forward sends 
coherent complexes to graded quasi-coherent complexes $M$ on $Z$ satisfying the following
two conditions.
\begin{enumerate}
\item For each $n$, the 
graded piece $M(n)$ is bounded coherent on $Z$.

\item For $n << 0$ we have $M(n)=0$.

\end{enumerate}

We let $\widehat{\Coh}(Z \times B\Gm) \subset 
\QCoh(Z \times B\Gm)$ be this full sub-$\s$-category. Note that it contains $\Coh(Z\times B\Gm)$,
the
graded coherent complexes on $Z$, but is strictly bigger as it also contains 
objects with infinitely many non-zero positive weight pieces. Symbolically, 
objects in $\Coh(Z\times B\Gm)$ will be denoted by 
$\sum_{i\in \ZZ}M(i).t^i$, and thus will be considered 
as Laurent formal power series with coefficients in $\Coh(Z)$. Another natural notation 
for $\widehat{\Coh}(Z\times B\Gm)$ is thus $\Coh(Z)[[t]][t^{-1}]$.

The pushfoward $\s$-functor 
$$p_* : \Coh([\VV^{(2)}(E[1])_Z/\Gm]) \to \Coh(Z)[[t]][t^{-1}]$$ then induces
a morphism on the corresponding Grothendieck groups
$$p_* : G_0^{\Gm}(\VV^{(2)}(E[1])_Z) \longrightarrow G_0(Z)[[t]][t^{-1}].$$
Finally, we can compose with the Grothendieck-Riemann-Roch transformation (see Appendix)
$$\tau_Z : G_0(Z)[[t]][t^{-1}] \longrightarrow H_{2*}^{BM}(Z,\Ql(-*))[[t]][t^{-1}],$$
where the right hand side is formal power series with coefficients in Borel-Moore homology. 
The composed
morphism will be denoted by 
$$\psi_Z : G_0^{\Gm}(\VV^{(2)}(E[1])_Z) \longrightarrow H_{2*}^{BM}(Z,\Ql(-*))[[t]][t^{-1}].$$

The following result is a version of the existence of Hilbert polynomials and is well known.

\begin{prop}\label{p2}
The morphism $$\psi_Z : G_0^{\Gm}(\VV^{(2)}(E[1])_Z) \longrightarrow 
H_{2*}^{BM}(Z,\Ql(-*))[[t]][t^{-1}]$$
factors through
the subgroup consisting of formal Laurent series 
$P(t) \in H_{2*}^{BM}(Z,\Ql(-*))[[t]][t^{-1}]$ such that 
there exist positive integers $n$ and $m$ with $(1-t^2)^n.t^m.P(t) \in H_{2*}^{BM}(Z,\Ql(-*))[t]$. In 
other words
$P(t)$ is always a rational function, possibly with poles at $t=-1, 0, 1$.
\end{prop}

\textit{Proof.} We start by assuming that $X$ has the resolution property: any coherent 
sheaf of $H^0(\OO_X)$-modules $\F$ is a quotient of a vector bundle $V \to \F$ on $X$. Under this
assumption, the perfect complex $E$ can be presented as the fiber of a morphism of
vector bundles $d : E_0 \to E_1$. In this case, we have a $\Gm$-equivariant closed embedding
$$j : \VV^{(2)}(E[1])_Z \hookrightarrow \VV^{(2)}(E_1)_Z,$$
setting over $Z$. The morphism $p_*$ can then be written as 
$q_*j_*$ where $q : \VV^{(2)}(E_1)_Z \to Z$ is the natural projection, and as $j_*$ preserves
bounded coherent complexes, the proposition reduces to the case where $E=E_1[-1]$ is a shift of 
a vector bundle. 
By homotopy invariance, the pull-back now induces an isomorphism
$$p^* : G_0^{\Gm}(Z) \simeq G_0^{\Gm}(\VV^{(2)}(E_1)_Z).$$
Moreover, by the projection formula (see Appendix) we have $p_*p^*(x)=p_*(1).x$, where 
we consider $G_0(Z)[[t]][t^{-1}]$ as a module over the ring
$K_0(Z)[[t]][t^{-1}]$ and $1$ is the class of the structure sheaf 
in $K_0^{\Gm}(\VV^{(2)}(E_1)_Z)$. 
Therefore, any element $p_*(y) \in G_0(Z)[[t]][t^{-1}]$ is of the form
$p_*(1).x$ with $x \in G_0(Z)[t,t^{-1}]$ in particular, 
this assertion signifies the existence of a positive integer $m$ such that 
$t^m.x$ is polynomial. Consequently, it suffices to
demonstrate the existence of another positive integer n such that
$(1-t^2)^n.Ch(p(1))$ is also polynomial. 
As the GRR transformation $\tau$ is a morphism of modules over the Chern character 
on $K$-groups, we see that it is enough to prove there is an integer 
$n$ such that $(1-t^2)^n.Ch(p_*(1))$ lies in 
$H^{2*}(Z,\Ql(*))[t]$, where $Ch : K_0(Z) \to H^{2*}(Z,\Ql(*))$ is the Chern character map to 
$\ell$-adic cohomology.

The element $p_*(1)$ is explicitly given by the formal power series 
$$S_{t^2}([E_1^\vee])=\sum_{i\geq 0}[Sym^i_{\OO_Z}(E_1^\vee)].t^{2i}$$
To prove that $(1-t^2)^n.Ch(p_*(1))$ is a polynomial for $n$ large enough, we can use
the splitting principle to reduce to the case where $E_1=\oplus_\alpha L_\alpha$ is a direct 
sum of line bundles.
In this case, we have $S_{t^2}([E_1^\vee])=\prod_{\alpha}S_{t^2}(L_\alpha)$, and we are reduced
to the case
where $E_1=L$ is a line bundle. Let $c=Ch(1-[L]) \in H^{2*}(Z,\Ql(*))$. The element $c$ is 
nilpotent, 
say $c^d=0$. We thus have
$$(1-t^2)^{d+1}.Ch(p_*(1))=(1-t^2)^{d+1}.(\sum_{i\geq 0}(1-c)^i.t^{2i}) = 
\sum_{i}(-1)^ic^i.(1-t^2)^{d-i}.t^i$$
which is indeed a polynomial.

Finally, for the general case, we can use the Chow envelope techniques of 
\cite{toengrr} in order
to produce a proper morphism $\pi : X' \to X$ such that, on one hand
$\pi_* : G_0(X'\times_X Z) \to G_0(Z)$ is surjective, and moreover $X'$ is smooth and 
possesses the
resolution property (in fact $X'$ can be taken to be a quotient stack
of a smooth quasi-projective variety by a finite group action).
\hfill $\Box$\\

\begin{rmk}
\emph{A slightly weaker conclusion of Proposition \ref{p2} almost holds in $G_0(Z)[[t]][t^{-1}]$, before
applying the GRR transformation to Borel-Moore homology. When $X$ is a derived algebraic
space, it holds as is. However, when $X$ is Deligne-Mumford the statement can not 
be true as stacky phenomena imply that the formal power series can also have poles
at roots of unity.}
\end{rmk}

One consequence of proposition \ref{p2} is that it is possible to evaluate
the series $\Eu_E(s)$ at $t=-1$, at least after removing the possible poles. 
For this, we remind the notion of critical value of a rational function $P \in K(t)$ with coefficients
in some field $K$.
By definition, the critical value of $P$ at $a \in K$ is the first non-zero Taylor coefficients
in a Laurent series expansion of $P$ at $a$. It is denoted by $P^*(a)$. In formula, it is given by 
$$P^*(a):=((t-a)^{-n}P(t))(a),$$
where $n$ is the valuation of $P$ at $a$ (by convention $P^*(a)=0$ if $P=0$). 
If $P$ is a polynomial, and $P(a)\neq 0$, then 
$P^*(a)=P(a)$, but when $a$ is a root of $P$ the critical value $P^*(a)$ is always non-zero
(except when $P=0$). One nice feature of critical values is that they define a morphism of multiplicative monoids
$K(t) \to K$, i.e. $(PQ)^*(a)=P^*(a).Q^*(a)$. However, the critical value map is of course 
not additive in general. \\

With these notions and definitions, we define Euler series, global Euler series and reduced Euler
numbers as follows.

\begin{df}\label{d6}
\begin{enumerate}
\item The \emph{Euler series of $s$}
is defined to be 
$$\Eu_E(s;t):=(-1)^d\psi(Td(E).Td(X)^{-1}\cap\tau_X(\Eu_E(s))) \in 
H_{2*}^{BM}(Z,\Ql(-*))[[t]][t^{-1}].$$ 
Here,  
$Td(E) \in H^{2*}(Z,\Ql(*))$ is the Todd class of the perfect complex $E$,
$Td(X)^{-1}$ is the inverse Todd class of $\T_X$ and $d$ is the rank of $E$.

\item When $Z=Z(s) \subset X$ is proper over $k$, the \emph{global Euler series of $s$} is defined to be 
$$Eu_E(s;t):=f_*(\Eu_E(s;t)) \in \Ql[[t]][t^{-1}],$$
where $f_ * :  H_{2*}^{BM}(Z,\Ql(-*)) \to \Ql$ is the pushforward in Borel-Moore homology along the structure map $f : Z \to 
\Spec\, k$.

\item When $Z$ is proper over $k$, the \emph{reduced Euler number of $s$} is defined to be
$$Eu^{r}_E(s):=Eu_E(s;t)^*(-1) \in \Ql.$$

\end{enumerate} 
\end{df}

\begin{rmk}
\begin{enumerate}
\item \emph{It is important to note that in point $(3)$ above, the number $Eu_E(s;t)^*(-1)$ 
is well defined
as Proposition \ref{p2} implies that $Eu_E(s;t) \in \Ql(t)$ and thus its critical value at $-1$ 
is defined.}

\item \emph{By construction $Eu_E(s;t)$ has coefficients in $\Ql$ but it is easy to 
see that it has rational coefficients. Moreover, when $X$ is a derived algebraic space instead
of a derived DM stack, $Eu_E(s;t)$ has in fact integral coefficients. However, 
in general, denominators appear in the push-forward $f_ * :  H_{2*}^{BM}(Z,\Ql(-*)) \to \Ql$, as
$f$ is not a representable morphism.}

\item \emph{The Todd classes in Definition \ref{d6} $(1)$ are only there in order to have a comparison with the 
usual Euler class of vector bundles via the Grothendiek-Riemann-Roch theorem (see below), 
but are not truly essential.}

\item \emph{Finally, we should also mention that it is also possible to define
a non-reduced Euler number by considering 
the value of $(1-t^2)^n.Eu_E(s;t)$ at $t=-1$, where $n$ is the least integer such that
$(1-t^2)^n.Eu_E(s;t) \in \Ql[t,t^{-1}]$.}

\end{enumerate}
\end{rmk}

\textbf{The bundle case.} Let us assume that $E$ is a vector bundle on $X$, and that $s=0$ 
to make things simpler. Then, $G_0^{\Gm}(\VV^{(2)}(E[1]) \simeq
G_0(X)[t,t^{-1}]$ as the projection $\VV^{(2)}(E[1]) \to X$ 
induces an isomorphism on the underlying truncated $DM$-stacks.
The refined
Euler class is explicitly given by the series
$$\Eu_E(0;t)=\sum_{i\in \ZZ}[Ext^i_{Sym_{\OO_X}(E^\vee)}(\OO_X,\omega_{E}^{-1})].t^i$$
This also be written a product of two Laurent polynomials
$$\Eu_E(0;t) = (\sum_{i} [H^i(\OO_X)].t^i)(\sum_{j}[\wedge_{H^0(\OO_X)}^{d-j}(H^0(E^\vee))].t^j),
$$
where $d$ is the rank of $E$. Theferore, setting $t=-1$ we get 
$$Eu_E(0;-1):=(-1)^d.\lambda_{-1}(E^\vee)\cap [X]_G^{vir},$$
where $\lambda_{-1}(E^\vee) = \sum_{i}(-1)^i[\wedge^iE^\vee] \in K_0(X)$ and 
$[X]_G^{vir}=\sum_i [H^i(\OO_X)] \in G_0(X)$
is the virtual class in $G$-theory of $X$. As the virtual class in homology is given by the so-called 
Kontsevich's formula (see \cite[Thm. 6.12]{khank})
$$[X]^{vir}=Td(X)^{-1}\cap \tau_{X}([X]^{vir}_G),$$
and as we have
$$Ch(\lambda_{-1}(E^\vee))=C_{top}(E).Td(E)^{-1},$$
we see that our Definition \ref{d6} gives that the reduced Euler number 
$Eu_E(0)$ equals $f_*(C_{top}(E)\cap [X]^{vir})$ for $f : X \to \Spec\, k$ the structure 
morphism, 
if this number is non-zero. In other words, the reduced Euler number coincide with the usual
Euler number of $E$ when it does not vanish. When the Euler number $f_*(C_{top}(E)\cap 
[X]^{vir})$ is zero, 
this means that $Eu_E(0;t)$ vanishes at $t=-1$, and the reduced Euler number is then equal to 
the
first non-zero derivative of the function $Eu_E(s;t)$ evaluated at $t=-1$. Note that, however, 
the non-reduced Euler number always coincides with $f_*(C_{top}(E)\cap [X]^{vir})$.

\subsection{Categories of shifted singularities}\label{seccatsing}

In this section we introduce a categorification of what we have seen so far. Namely
we introduce certain categories of singularities for shifted potentials, that is
for functions of cohomological degree $-1$. In the same manner that 
categories of singularities for a usual potential $f$ reflects the singularities of $f$, 
these categories of shifted singularities reflects singularities of shifted functions
in the sense that they see the locus where the shifted function is not quasi-smooth.
The results contained in this section are independant of the rest of the paper and will
not be used in the forthcoming sections.

Let $Y$ be a quasi-smooth derived Deligne-Mumford stack and $f : Y \to \mathbb{A}^1[-1]:=\Spec\,
k[\epsilon_{-1}]$ be a function of degree $-1$ (here $\epsilon_{-1}$ sits is cohomological
degree $-1$). The function $f$ can also be considered as a morphism of $\OO_Y$-algebras
$$\OO_Y[\epsilon_{-1}] \to \OO_Y,$$
that sends $\epsilon_{-1}$ to $f$. Therefore, any perfect complex $M$ on $Y$
is also a coherent $\OO_Y[\epsilon_{-1}]$-module. We will say that $M \in \Parf(Y)$ 
is $f$-perfect if the complex of sheaves
$$\rch_{\OO_Y[\epsilon_{-1}]}(\OO_Y,E)$$
is locally cohomologically bounded on $Y$. We denote
by $\Parf(Y,f) \subset \Parf(Y)$ the full sub-$\s$-category of $f$-perfect complexes.

We then set
$$\Sing(Y,f):=\Parf(Y)/\Parf(Y,f),$$
where the quotient is taken in $k$-linear dg-categories up to Morita equivalences.
In the same manner, when $Y$ comes equipped with a $\Gm$-action, and $f : Y \to \mathbb{A}^1[-1]$
is $\Gm$-equivariant, for the standard $\Gm$-action on  $\mathbb{A}^1[-1]$ of weight $1$, 
we have an equivariant version
$$\Sing^{\Gm}(Y,f):=\Parf^{\Gm}(Y)/\Parf^{\Gm}(Y,f).$$
We note that $\Sing(Y,f)$ is naturally enriched over $\Parf(Y)$, and thus can be
considered as a quasi-coherent stack of dg-categories over $Y$ (see 
\cite{MR2957304} for the general notion of quasi-coherent 
stacks of dg-categories). As such, 
it is important to notice the following fact. \\

\begin{prop}\label{p3}
As a quasi-coherent stack of dg-categories on $Y$, $\Sing(Y,f)$ (resp. $\Sing^{\Gm}(Y,f)$
in the equivariant case) is supported
on the (closed) locus where $f : Y \to \mathbb{A}^1[-1]$ is not quasi-smooth.
\end{prop}

\textit{Proof.} We have to show that if $f : Y \to \mathbb{A}^1[-1]$ is quasi-smooth
then $\Sing(Y,f)\simeq 0$, or equivalently that all perfect complexes on $Y$
are $f$-perfect. As $f$ is quasi-smooth $\OO_Y$, locally on $Y$, admits a finite
resolution by free $k[\epsilon]_-{-1}$-dg-modules. 
This implies that 
for any perfect complex $E$ on $Y$, locally on $Y$, $E$ admits a finite resolution
by free $k[\epsilon_{-1}]$-dg-modules. 
As $\rh_{k[\epsilon_{-1}]}(k,k[\epsilon]) \simeq k$, 
the complex of sheaves $\rch_{\OO_Y[\epsilon_{-1}]}(\OO_Y,E)$
is indeed locally cohomologically bounded. This is precisely the statement that $\Parf(Y,f)=\Parf(Y)$.
\hfill $\Box$ \\

In our specific situation we apply this to $Y=\VV(E[1])$ and 
$f_s : \VV(E[1]) \to \mathbb{A}^1[-1]$
is the function induced from the section $s$. Indeed, $s : \OO_X \to E$ induces
a morphism $s[1] : \OO_X[1] \to E[1]$, and thus a morphism on the corresponding 
linear stacks
$$\VV(E[1]) \to \VV(\OO_X[1]) \simeq X \times \mathbb{A}^1[-1],$$
whose second projection is, by definition, our function $f_s$. This function is moreover compatible with 
the $\Gm$-action by weight $1$. We thus have a dg-category
$$\Sing^{\Gm}(\VV(E[1]),f_s) = \Parf^{\Gm}(\VV(E[1]))/\Parf^{\Gm}(\VV(E[1]),f_s).$$

We want to see the dg-category $\Sing^{\Gm}(\VV(E[1]),f_s)$, together with its natural object
$\omega_E^{-1} \in \Sing^{\Gm}(\VV(E[1]),f_s)$, 
as some sort of categorification of the pair
$(G_0^{\Gm}(\VV^{(2)}(E[1])_Z),\Eu_E(s)) $ where $\Eu_E(s)$
is our refined Euler class.
Indeed, we will exhibit below a natural $\s$-functor
$$|-|^t : \Sing^{\Gm}(\VV(E[1]),f_s) \longrightarrow 
\Coh_{\heartsuit}^{\Gm}(\VV^{(2)}(E[1]))_Z,$$
where $\Coh^{\Gm}_{\heartsuit}(\VV^{(2)}(E[1]))_Z$ is the abelian category of $\Gm$-equivariant
sheaves on $\VV^{(2)}(E[1])$ set theoretically supported on $Z$, such that $[|\omega_E^{-1}|^t] = \Eu_E(s)$
in $G_0^{\Gm}(\VV^{(2)}(E[1])_Z)$. \\

\noindent \textbf{Construction of the functor $|-|^t$.} In order to define the functor $|-|^t$, first we notice that we have to distinguished two objects in
$\OO_{\VV(E[1])}[\epsilon_{-1}]$-modules. We have 
$\OO_{\VV(E[1])}$, endowed with its $\OO_{\VV(E[1])}[\epsilon_{-1}]$-module structure
coming from the augmentation $\epsilon_{-1}\to 0$. We also have
$\OO_{\VV(E[1])}$ with its natural $\OO_X$-algebra morphism
$\OO_{\VV(E[1])}[\epsilon_{-1}] \to \OO_{\VV(E[1])}$ sending $\epsilon_{-1}$
to $f_s$. These are two different $\OO_{\VV(E[1])}[\epsilon_{-1}]$-module structures
on the same object $\OO_{\VV(E[1])}$-module $\OO_{\VV(E[1])}$ itself. 
More generally, 
for any $M \in \Parf^{\Gm}(\VV(E[1]))$ we have two different 
$\OO_{\VV(E[1])}[\epsilon_{-1}]$-module structures on $M$, denoted by 
$M\{0\}$ and $M\{f_s\}$. 
For $M\{0\}$, $\epsilon_{-1}$ acts by $0$, whereas
for $M\{f_s\}$ the element $\epsilon_{-1}$ acts by multiplication by $f_s$.
For the sake of simplicity, $M\{0\}$ will simply be denoted by $M$, but we will
keep the notation $M\{f_s\}$ for the second one.

We start by constructing a dg-functor
$$|-| : \Parf^{\Gm}(\VV(E[1])) \longrightarrow \Parf^{gr}(Sym_{\OO_X}(E[1])[u]),$$
where $u$ is a variable in cohomological degree $2$ and weight $-1$ (note here
that $E$ has weight $1$), and 
$\Parf^{gr}(Sym_{\OO_X}(E[1])[u])$ are sheaves
of perfect graded $Sym_{\OO_X}(E[1])[u]$-dg-modules.
This dg-functor is essentially given by negative cyclic homology, and 
the weight $n$ piece of $|M|$ is by definition 
$$|M|(n):=\underline{Hom}_{\OO_{\VV(E[1])}[\epsilon_{-1}]}(\OO_{\VV(E[1])}(n),M\{f_s\}),$$
where $\OO_{\VV(E[1])}(n)$ is $\OO_{\VV(E[1])}$ with weight shifted by $n$.
The action of $u$ is given by the usual identification
$$\oplus_n |\OO_{\VV(E[1])}(n)| \simeq \OO_{\VV(E[1])}[u].$$
By definition, the full sub-dg-category $\Parf^{\Gm}(\VV(E[1]),f_s)$ is sent to (locally) cohomologically bounded
objects, and thus its image consists of $Sym_{\OO_X}(E[1])[u]$-dg-modules
$E$ for which $colim_{i}E[2i] \simeq 0$, where the transition morphisms
in the colimit $E[2i] \to E[2i+2]$ are given by mutiplication by $u$. This means that
$|-|$, composed with the canonical localization dg-functor
$$\Parf^{gr}(Sym_{\OO_X}(E[1])[u]) \to \Parf^{gr}(Sym_{\OO_X}(E[1])[u,u^{-1}])$$
sends $\Parf^{\Gm}(\VV(E[1]),f_s)$ to $0$.
We thus obtain a well defined dg-functor
$$|-| : \Sing^{\Gm}(\VV(E[1]),f_s) \longrightarrow \Parf^{gr}(Sym_{\OO_X}(E[1])[u,u^{-1}]).$$
Now, we recall that there exists a natural tensor autoequivalence of 
graded complexe, called \emph{red shift}, which sends a complex $E$ pure
of weight $n$ to $E[2n]$ with the same weight $n$. This autoequivalence, 
implies
that $\Parf^{gr}(Sym_{\OO_X}(E[1])[u,u^{-1}])$ is equivalent to \\
$\Parf^{gr}(Sym_{\OO_X}(E[-1])[t,t^{-1}])$ where now $t$ sits in cohomological degree $0$ and
weight $-1$. This last dg-category is naturally equivalent to $\Parf(Sym_{\OO_X}(E[-1]))$, 
non-graded dg-modules over $Sym_{\OO_X}(E[-1])$. Clearly, the composed
dg-functor 
$$|-| : \Sing^{\Gm}(\VV(E[1]),f_s) \longrightarrow \Parf(Sym_{\OO_X}(E[-1]))$$ takes
it values in the full sub-dg-category $\Parf(Sym_{\OO_X}(E[-1]))_Z$ of objects supported
on $Z=Z(s) \subset X$.
We can then compose with the total cohomology functor
$$H^* : \Parf(Sym_{\OO_X}(E[-1]))_Z \longrightarrow 
\Coh^{\Gm}_{\heartsuit}(\VV^{(2)}(E[1]))_Z,$$
that sends $M$ to $H^*(M)$ considered as graded $Sym_{\OO_X}(H^1(\T_X))$-module. We get this
way the desired $\s$-functor
$$|-|^t :  \Sing^{\Gm}(\VV(E[1]),f_s) \longrightarrow \Coh_{\heartsuit}^{\Gm}(\VV^{(2)}(E[1]))_Z.$$

\bigskip

Note that $|-|^t$ is not a stable functor nor a dg-functor, and thus does not preserve exact 
triangles. In particular, this $\s$-functor does not induce a
well defined map on $K$-groups.
We consider the determinant line bundle $\omega_E^{-1}=det(E)^{-1}$ on $X$, and its
inverse image on $\VV(E[1])$ which sits in $\Parf^{\Gm}(\VV(E[1]))$. The image
of this object by the quotient map $\Parf^{\Gm}(\VV(E[1])) \to \Sing^{\Gm}(\VV(E[1]))$
will still be denoted by $\omega_E^{-1}$. \\

\begin{prop}\label{p4}
With the above notations, we have 
$$[|\omega_E^{-1}|^t] = \Eu_E(s)$$
in $G_0^{\Gm}(\VV^{(2)}(E[1])_Z)$.
\end{prop}

\textit{Proof.} This follows from the explicit description of the dg-functor
$$|-| : \Parf^{\Gm}(\VV(E[1])) \longrightarrow \Parf^{gr}(Sym_{\OO_X}(E[1])[u]).$$
An object $M \in \Parf^{\Gm}(\VV(E[1]))$ is given by a graded $Sym_{\OO_X}(E[1])$-module.
The dg-functor $|-|$ sends such an $M$ to the graded $Sym_{\OO_X}(E[1])[u]$ whose
weight $p$ pieces are given by
$$|M|(p)=\prod_{i\geq p}M(i)[2p-2i],$$
where the differential is given by the total differential, sum of 
the cohomological differential and multiplication by the degree $-1$ function $f_s$. 
The action of $u$ is itself given by the canonical embedding of sub-products
$$|M|(p)=\prod_{i\geq p}M(i)[2p-2i] \hookrightarrow \prod_{i\geq p-1}M(i)[2p-2i]=|M|(p-1)[2].$$
In particular, $|M|^t$ is simply given by the Tate realization (see \cite{cptvv})
of the graded mixed object $M\{f_s\}$
$$|M|^t := colim_p H^*(\prod_{i\geq p}M(i)[-2i]),$$
considered now as a graded-$H^*(Sym_{\OO_X}(E[-1]))$-module.
When applied to $M=\omega_{E}^{-1}$, 
we get that $|\omega_E^{-1}|^t$ is represented by
$Sym_{\OO_X}(E[-1]) \otimes_{\OO_X}\omega_{E}^{-1}$
where the differential on $Sym_{\OO_X}(E[-1])$ is the sum of the cohomological 
differential and the multiplication by $f_s$. This clearly 
is $\underline{Hom}_{\OO_{\VV(E^\vee)}}(e_*\OO_X,s_*\omega_E^{-1})$
as a module over $$Sym_{\OO_X}(E[-1])\simeq 
\underline{End}_{\OO_{\VV(E^\vee)}}(e_*\OO_X),$$ hence the statement.
\hfill $\Box$ \\

\begin{rmk}
\begin{enumerate}
\item \emph{The above Proposition \ref{p4} should be related to the fact that 
the usual virtual structure sheaf of $Z(s)$, when $E$ is a bundle, can be 
recovered from the theory of vanishing cycles of the function $\VV(E^\vee) \to \mathbb{A}^1$
induced by $s$ (see for instance 
\cite[Appendix A]{MR3667216})). Here, the role of vanishing cycles has been replaced by the
category of shifted singularities, which is relevant by the well known relations
between vanishing cycles and categories of singularities (\cite{brtv}). However, we do not try here
to establish a precise relation between $\Sing^{\Gm}(\VV(E[1]),f_s)$ and the
vanishing cycles of $\VV(E^\vee) \to \mathbb{A}^1$, even though we think that 
some precise statement could be written down. }

\item \emph{Proposition \ref{p4} is not only another way to understand
the class $\Eu_E(s)$, it also provides a nice representative of this class
as a nice object $\OO_{\VV(E[1])} \in \Sing^{\Gm}(\VV(E[1]),f_s)$, 
the image of the structure sheaf of $\VV(E[1])$ by the canonical projection. One advantage 
of this object is that it has nice functorial properties in terms of $X$, $E$ and $s$. Such properties
are not shared by the corresponding sheaf of graded $Sym_{\OO_X}(H^1(\T_X))$-modules 
as this arises as the total cohomology and thus will not be stable under various base changes.}

\item \emph{Finally, our $\Sing(Y,f)$ is also closely related to the localisation 
of the virtual class along a cosection of \cite{KL}. Indeed, the non-equivariant 
version of our $\s$-functor $|-|^t$ is an $\s$-functor
$$|-|^t : \Sing(Y,f) \longrightarrow \Coh(Y)_Z$$
where $Z$ is the zero locus of $f$ in $Y$. This $\s$-functor $|-|^t$ simply sends
$M$ to its periodic cyclic homology $\underline{Hom}_{\OO_Y[\epsilon_{-1}]}(\OO_Y,M{f})[u^{-1}]$,
which is a $2$-periodic complex with coherent cohomologies on $Y$, supported on $Z$.
We thus have a canonical element 
$$[|M|^t]=[H^0(\underline{Hom}_{\OO_Y[\epsilon_{-1}]}(\OO_Y,M{f})[u^{-1}])]-
[H^1(\underline{Hom}_{\OO_Y[\epsilon_{-1}]}(\OO_Y,M{f})[u^{-1}]))] $$
inside $G_0(Z)$.
Applied to $M=\OO_Y$ the resulting class $[|\OO_Y|^t]$ seems very 
closely related to the localized
virtual class of \cite{KL} (though we do not investigate a precise formula in this work).}

\end{enumerate}
\end{rmk}

\section{Applications}

In this final section, we will apply the results and constructions of the previous section to the specific 
case of the derived stack of stable maps endowed with an induced vector field as in Definition \ref{d3}. 

\subsection{Categories associated to a vector field}

Let $X$ be a smooth and proper algebraic variety, and $\nu \in  H^0(X,\T_X)$ a global vector field on $X$. We will
construct three invariants of $\nu$, based on the induced vector field $\nu_{0,2}$ on $\RMgn_{0,2}(X)$. These
are of three different levels, 2-categorical, categorical and numerical/homological, and are
somehow related by some form of decategorification. They are all incarnations of a possible algebraic
version of Morse homology, or more precisely the specific point of view based
on \emph{moduli broken lines} exposed in \cite{lurie2018associative}. We note also that the
monoid $M_{0,2}(X)$ has been first studied in \cite{MR2361096}, and that our constructions in this section
might be related to the path Lie algebroid introduced in \cite{MR2361096}.

We first set $\M:=\coprod_{\beta} \RMgn_{0,2}(X,\beta)$, where $\beta$ runs over all 
curve classes in $H_2(\overline{X},\Ql(-1))$. The additive monoid of curve classes will be denoted by $N \subset H_2(\overline{X},\Ql(-1))$. 
As a side comment, $\M$ is not quasi-compact, and when considering non-quasi-compact derived stacks
we will only consider coherent complexes supported on quasi-compact substacks. In particular, 
our notation for $\Coh(\M)$ will mean $\oplus_\beta \Coh(\M)$ rather than the
infinite product $\prod_\beta \Coh(\M)$. The same convention will hold each time we
consider $\Coh$ of some derived stack which is an infinite disjoint union of 
quasi-compact derived stacks.

The two marked points of prestable curves will be denoted by $0$ and $\s$, and we think of a stable
map $f : C \to X$ as a path from $f(0)$ to $f(\s)$. The gluing construction provides a composition morphism
$$\mu : \M \times_X \M \longrightarrow \M,$$
where the two maps $\M \to X$ defining the fibered product on the left hand side
are the two evaluation morphisms. This composition makes $\M$ into an associative monoid (without units)
in the $\s$-category $\dSt/(X\times X)$ of derived stacks over $X\times X$ endowed with 
the monoidal structure given by convolution. Another possible terminology here could be that 
$\M$ is a categorical object in derived stacks whose object of objects is $X$. As such, 
$\M$ is very similar to the path groupoid of a topological space.

The vector field $\nu$ induces a vector field $\nu_{0,2}$ on $\M$ as explained in Proposition \ref{p1}. Moreover, $\nu$ and
$\nu_{0,2}$ are compatible with the evaluation $\M \to X \times X$ as well as with the composition 
morphism $\mu$ in an obvious sense. As a result, the derived scheme of zeros $\nu_{0,2}$, denoted
by $Z \subset \M$, comes equipped with an induced composition map
$$\mu : Z \times_{Z_0} Z \longrightarrow Z$$
making it again into an associative monoid (without units) in $\dSt/(Z_0\times Z_0)$, 
where $Z_0 \subset X$ is the derived subscheme of zeros of $\nu$. \\

\textbf{The 2-category of algebraic trajectories.}
Let $ev : Z \to Z_0 \times Z_0$ be the evaluation morphism. We consider
$\Coh(Z/ev) \subset \Coh(Z)$, the full sub-dg-category 
of coherent complexes which are of finite Tor-dimension over $Z_0 \times Z_0$
(i.e. that are relatively perfect with respect to $ev$). This is a $\Parf(Z_0 \times Z_0)$
linear dg-category, and will be considered as a quasi-coherent sheaf of dg-categories
over $Z_0\times Z_0$ (see \cite{MR2957304}). We note that the monoid structure on $Z$ induces a monoid structure on $\Coh(Z/ev)$, 
as an object in quasi-coherent sheaves of dg-categories over $Z_0 \times Z_0$ endowed with 
the convolution monoidal structure. In concrete terms, this monoid structure is given by 
$$\Coh(Z/ev) \otimes_k \Coh(Z/ev) \longrightarrow \Coh(Z/ev)$$
sending an object $(E,E')$ to $\mu_*(d^*(E \boxtimes E'))$, where 
$d : Z\times_{Z_0} Z \to Z \times Z$ is the natural embedding. Note that $d^*$
does preserve coherent complexes which are relatively perfect over $Z_0 \times Z_0$ even though
it does not preserve coherent complexes in general.
The dg-functor $\Coh(Z/ev) \otimes_k \Coh(Z/ev) \longrightarrow \Coh(Z/ev)$ is linear over the
multiplication (tensor product) $\Parf(Z_0) \otimes \Parf(Z_0) \to \Parf(Z_0)$ and thus
descends to the mentioned monoid structure
$$\Coh(Z/ev) \otimes_{\Parf(Z_0)}\Coh(Z/ev) \to \Coh(Z/ev).$$

We therefore have defined a "linear $2$-category" (without units), whose "set" of object is $Z_0$,  
whose dg-category of morphisms is given by $\Coh(Z/ev)$. 

\begin{df}\label{d7}
The \emph{$2$-category of algebraic trajectories of $\nu$ on $X$} is
$\Coh(Z/ev)$, considered as a monoid in $Dg^c(Z_0 \times Z_0)$. It is denoted by 
$\Tra(X,\nu)$.
\end{df}

A very important special case is when $Z_0$ is of dimension zero, and thus is of the form
$\coprod_{i \in I} \Spec\, A_i$ where the $A_i$'s are local artinian $k$-dg-algebras $A_i$, and $I$ is the set of 
zeros of $\nu$ in $X$.
The derived stack $Z$ then splits as $\coprod_{(i,j)}Z_{i,j}$, where $Z_{i,j}$ lives over $\Spec\, A_i \times
\Spec\, A_j$. The $2$-category $\Tra(X,\nu)$ is then given by a concrete $2$-category whose set of objects
is $I$. The dg-category of morphisms from $i$ to $j$ is given by
$$\Tra(X,\nu)(i,j)=\Coh(Z_{i,j}/ev),$$
the dg-category of coherent complexes on $Z_{i,j}$ which are relatively perfect over $A_i \otimes A_j$.
Composition of morphisms is given by the gluing maps
$$\Coh(Z_{i,j}/ev) \otimes_{A_j} \Coh(Z_{j,k}/ev) \longrightarrow \Coh(Z_{i,k}/ev).$$
By our convention, $Z$ is itself the union of $Z(\beta)$, where $\beta \in N$ runs over curve classes. This
decomposition provides an $N$-graduation on the dg-categories of morphisms $\Coh(Z_{i,j}/ev)$, 
in a way that the composition above is $N$-graded. To be more precise, we have
$\Coh(Z_{i,j}/ev) \simeq \oplus_{\beta \in N}\Coh(Z_{i,j}(\beta)/ev)$, and
the composition sends the component $$\Coh(Z_{i,j}(\beta)/ev) \otimes_{A_j} \Coh(Z_{j,k}(\beta')/ev)$$
to $\Coh(Z_{i,k}(\beta+\beta')/ev)$.

Finally, when $Z_0$ is moreover reduced and each zero is rational over $k$,
we have $A_i\simeq k$ for all $i$, and
$\Tra(X,\nu)$ becomes a genuine $k$-linear 2-category (but still without units).\\

\textbf{The category of algebraic trajectories.}
We now construct a (dg) 1-categorical analogue of the $2$-category $\Tra(X,\nu)$. It can probably be obtained by
considering the $\ell$-adic realization of $\Tra(X,\nu)$ by applying the techniques of \cite{brtv} levelwise
on dg-categories of morphisms (or rather on $\Coh(Z/ev)$ as a dg-category over $Z_0 \times Z_0$).
We will however construct it directly using the $\ell$-adic formalism. 
We come back to the evaluation map $ev : Z \to Z_0 \times Z_0$. We consider
$\mathsf{Traj}(X,\nu):=ev_*ev^!(\Ql)$, the relative Borel-Moore homology of $Z$ over $Z_0 \times Z_0$. This
is a constructible $\ell$-adic complex on $Z_0 \times Z_0$. As $ev$ is proper, the fiber 
of $\mathsf{Traj}(X,\nu)$ at a geometric point $(x,y) : \Spec\, K \in Z_0\times Z_0$ is 
the Borel-Moore homology complex of $Z_{(x,y)}:=Z \times_{Z_0\times Z_0} \Spec\, K$.

As in the previous $2$-categorical setting, we endow the $\s$-category  $\D(Z_0\times Z_0,\Ql)$ of
constructible $\ell$-adic complexes with the natural convolution tensor product (using $*$-pullbacks). 
We claim that $\mathsf{Traj}(X,\nu)$ is then endowed with a natural structure of an associative (non-unital)
monoid in $\D(Z_0\times Z_0,\Ql)$. The multiplication map for this monoid structure is
induced by the composition map  $\mu : Z\times_{Z_0} Z \to Z$ and its natural direct image in
relative homology
$$\mu_* : p_*(p!(\Ql)) \to ev_*ev^!(\Ql)$$
where $p : Z\times_{Z_0} Z \to Z_0 \times Z_0$ is the natural structure map (so that
$\mu$ is a morphism over $Z_0\times Z_0$).

We get this way a monoid in $\D(Z_0\times Z_0,\Ql)$, without units, 
whose underlying object is $\mathsf{Traj}(X,\nu)$. By definition, this is a $\Ql$-linear dg-category whose
"set of objects" is $Z_0$. 

\begin{df}\label{d8}
The \emph{dg-category of algebraic trajectories of $\nu$ on $X$} is
$\mathsf{Traj}(X,\nu)$, considered as a monoid in $\D(Z_0 \times Z_0,\Ql)$. 
\end{df}

Again, as in the $2$-categorical case, when $Z_0$ is of the form $\coprod_\Spec\, k$, 
$\mathsf{Traj}(X,\nu)$ is a genuine dg-category whose set of objects is $Z_0$, 
(possibly endowed with a Galois action when $k$ is not algebraically closed). \\

\textbf{Numerical and homological invariants.} We now go one step further in the decategorification and 
define numerical and homological invariants. For each curve class $\beta \in N$ we have
$ev_\beta : Z_\beta \to Z_0 \times Z_0$, the restriction of the evaluation map to the $\beta$
component. The direct image by $ev_{\beta}$ of the Euler series $\Eu_E(s,t) \in H_{2*}(Z_\beta,\Ql(-*))((t))$
of Definition \ref{d6}, where $E=\T_X$ and $s=\nu_{0,2}$, provides a natural element 
$$\widetilde{d_\beta} \in (H_{2*}(Z_0,\Ql(-*))\otimes H_{2*}(Z_0,\Ql(-*)))((t)).$$
This element should be thought as some form of Morse differential and as being the algebraic analogue
of the usual differential in Morse homology (see \cite{bot}). For instance,
when $Z_0$ is smooth, $\widetilde{d_\beta}$ can be considered, via Poincaré duality,
as a $\Ql((t))$-linear endomorphism of $H_{2*}(Z_0,\Ql(-*))((t))$. It is not true in general
that $\widetilde{d_\beta}^2=0$ so this endomorphism is not quite a differential. 

The situation seems slightly better when $X$ is convex, so that $\RMgn_{0,2}(X)$ is smooth, and
the Laurent series $\widetilde{d_\beta}$ is in fact a polynomial in $t$. In this case, we can define an endomorphism
by specializing $t=-1$
$$d_\beta := \widetilde{d_\beta}(t=-1) \in H_{2*}(Z_0,\Ql(-*))\otimes H_{2*}(Z_0,\Ql(-*)).$$
This element is the direct image by the evaluation map 
$$ev : Z_\beta \times_{Z_0} Z_\beta \longrightarrow Z_0 \times Z_0,$$
of the localized top Chern class of $\RMgn_{0,2}(X,\beta)$ along the section $\nu_{0,2}$. In particular, 
this element is of homological degree $0$. When $Z_0$ is smooth and $0$-dimensional (i.e. $Z_0$ 
is smooth as a derived scheme of zeros of $\nu$), so $\nu$ has
only simple and isolated zeros on $X$, we get this way an endomorphism of the vector space
spanned by the set of zeros $I$
$$d_\beta : I\otimes\Ql \longrightarrow I\otimes \Ql.$$
To be more precise, we can set $I$ to be the zero of geometric zeros (i.e. geometric points
of $Z_0$), so that it comes equipped with an action of the absolute Galois group $G$ of $k$. The endomorphism
$d_\beta$ is then $G$-equivariant.

Again, $d_\beta$ is not quite a differential, but there is an explicit formula for its square, 
very closely related to the splitting axiom in GW theory of \cite{Manin}. To get this formula, we introduce
another derived Deligne-Mumford stack $W_\beta$, together with a 
natural morphism $W_\beta \to Z_0 \times Z_0$. It is defined as follows. 
We let $\M_\beta:=\RMgn_{0,2}(X,\beta)$, and consider the projection to the stack 
of prestable curves of genus zero with two marked points $\M_{0,2}^{pre}$.
$$p : \M_\beta \to \Mgn_{0,2}^{pre}.$$
On the level of prestable curves we 
have the gluing morphism 
$$\mu : \left(\Mgn_{0,1}^{pre} \times \Mgn_{0,3}^{pre} \right) \coprod 
\left( \Mgn_{0,2}^{pre} \times \Mgn_{0,2}^{pre} \right) \coprod 
\left( \Mgn_{0,3}^{pre} \times \Mgn_{0,1}^{pre}\to \Mgn_{0,2}^{pre} \right).$$
The left hand side will be denoted by $\Mgn_{0,2}^{\s}$, and will be called
the moduli of \emph{decomposed prestable curves} (of genus $0$ with
two marked points). We then
form the cartesian square
$$\xymatrix{
\M_\beta^{\s}:= \M_\beta \times_{\Mgn_{0,2}^{pre}} \times \Mgn_{0,2}^{\s} \ar[r]
\ar[d] & \M_\beta \ar[d] \\
\Mgn_{0,2}^{\s} \ar[r] & \Mgn_{0,2}^{pre}.
}$$

The derived stack $\M_\beta^{\s}$ is by definition the boundary at $\s$ of $\M_\beta$, also 
called the \emph{derived stack of decomposed stable maps}. Note that $\mu$ is
a finite, lci and unramified morphism, so that $\M_\beta^{\s}$ is again a quasi-smooth proper derived Deligne-Mumford
stack. Its image in $\M_\beta$ is the canonical divisor of $\M_\beta$.
The vector field $\nu_{0,2}$ being vertical with respect to the projection $p$, it induces 
a natural vector field $\nu_{0,2}\oplus 0$ on $\M_\beta^{\s}$. Its derived stack of zeros
is denoted by $W_\beta \subset \M_{\beta}^{\s}$, and fits in a cartesian square
$$\xymatrix{
W_\beta  \ar[r]
\ar[d] & Z_\beta \ar[d] \\
\T\Mgn_{0,2}^{\s}[-1] \ar[r] & \T\Mgn_{0,2}^{pre}[-1],
}$$
where $\T Y[-1]=\Spec\, Sym(\LL_{Y}[1])$ is the derived stack
of zeros of the zero vector field on a given derived stack $Y$. 
Note that, as $X$ is convex, the morphism $p : \M_\beta \to \Mgn_{0,2}^{pre}$
is smooth. This implies that the induced morphism on zeros of vector fields $Z_\beta \to \T\Mgn_{0,2}^{pre}[-1]$
is itself quasi-smooth. As a result, $W_\beta$ is a quasi-smooth and closed substack of $\M_\beta$. It 
comes equipped with the induced projection $ev : W_\beta \to Z_0 \times Z_0$. The direct
image of the virtual class of $W_\beta$ by $ev$ produces another endomorphism of $I\otimes\Ql$ denoted
by
$$d_\beta^{\s}:=ev_*([W_\beta]^{vir}) : I\otimes\Ql \longrightarrow I\otimes \Ql.$$
Finally, we will also need to define two other objects coming from stable maps with $1$ and $3$ marked points. 
For $\beta \in N$, we have
$\RMgn_{0,1}(X,\beta)$, together with its evaluation map $\RMgn_{0,1}(X,\beta) \to X$. 
This induces a morphism  on the level of derived stacks of zeros of $\nu_{0,1}$
$$\mathcal{E}_\beta \to Z_0,$$
which is a proper and quasi-smooth morphism. The direct image of the virtual fundamental 
class $[\mathcal{E}_\beta]^{vir} \in H_0(\mathcal{E}_\beta,\Ql)$ defines
an element 
$$e_\beta \in H_0(Z_0,\Ql)\simeq I\otimes \Ql.$$
In the same manner, we consider $\RMgn_{0,3}(X,\beta)$ with its evaluation to 
$X^3$. On the level of derived stacks of zeros of $\nu_{0,3}$, we get 
another proper and quasi-smooth morphism
$$\mathcal{F}_\beta \to Z_0^3.$$
The direct image of the virtual class by the above morphism produces an element
$f_\beta \in (I\otimes \Ql)^3$. By Poincaré duality, this element can be paired with
an element $x \in I\otimes \Ql$ in two different manners, $x.f_\beta$ and $f_\beta.x$, in order
to get an endomorphism of $I\otimes \Ql$.

The following proposition is the analogue of the splitting/associativity axiom in GW theory (see e.g. \cite[Def. 7.1-(3)]{Manin}, taking into account
our vector field $\nu$.

\begin{prop}\label{p5}
Assume, as before, that $X$ is convex and that the zeros of $\nu$ are isolated and simple. Then,
for each $\beta \in N$
we have
$$\sum_{\beta_1 + \beta_2=\beta}d_{\beta_1}.d_{\beta_2} + e_{\beta_1}.f_{\beta_2}+f_{\beta_1}.e_{\beta_2}
= d_{\beta}^{\s}.$$
\end{prop}

\textit{Proof.} This follows easily from the natural identification of $\M_\beta^{\s}$ with 
$$\coprod_{\beta_1 + \beta_2=\beta} \left( \M_{\beta_1} \times_X \M_{\beta_2} \right) \coprod 
\left( \RMgn_{0,1}(X,\beta_1)\times_X \RMgn_{0,3}(X,\beta_2) \right) $$
$$\coprod  \left(
\RMgn_{0,3}(X,\beta_1)\times_X \RMgn_{0,1}(X,\beta_2) \right),$$
and the fact that the identification respects the induced vector fields.
\hfill $\Box$ \\

\noindent \textbf{An easy example.} We finish this section by the most simple example, namely $X=\mathbb{P}^1$ and $\beta=1$. In this case, 
$\Mgn_{0,2}(X,\beta) \simeq \mathbb{P}^1 \times \mathbb{P}^1$, the isomorphism being given by 
the evaluation map. Given $(a,b) \in \mathbb{P}^1 \times \mathbb{P}^1$ outside of the diagonal,
the corresponding stable map $C \to \mathbb{P}^1$ is given by any linear endomorphism of $\mathbb{P}^1$
sending $0$ to $a$ and $\s$ to $b$. When $a=b$, the corresponding stable map is such that 
$C$ has two components, one with the two marked points entirely contracted sent to $a$, and a second
one identified with $\mathbb{P}^1$ with the identity map. Let $\nu = z\partial_z$ be the standard
vector field on $\mathbb{P}^1$ with two simple zeros at $0$ and $\s$. The induced vector field
$\nu_{0,2}$ is the direct sum $\nu \oplus \nu$ as a vector field on $\mathbb{P}^1 \times \mathbb{P}^1$.
This vector field has $4$ simple zeros, namely $(0,0)$, $(0,\s)$, $(\s,0)$ and $(\s,\s)$. As a result, 
$I\otimes\Ql=\Ql.0 \oplus \Ql.\s$, and the endomorphism $d_1$ defined above
is given by 
$$d_1(0)=d_1(\s)=0+\s.$$
Therefore, we see that the equation satisfied by $d_1$ is here $d_1^2=2.d_1$.

\subsection{An algebraic dynamical Zeta function}\label{sectionzeta}

We now turn to the case of genus $1$. Let $\nu$ be again a vector field on 
a smooth and proper scheme $X$. For each curve class  $\beta \in N$ we 
have the total Euler series of the vector field $\nu_{1,0}$ 
on $\RMgn_{1,0}(X,\beta)$ of stable maps of genus $1$ with no marked points. 
The zeros of $\nu_{1,0}$ corresponds to "periodic orbits" of $\nu$ in $X$, 
that is elliptic curves $f : E \to X$ whose image is globally invariant by 
$\nu$. Let us fix a polarization $\omega$ on $X$, so that 
we have a degree map $\beta \mapsto |\beta|$
$$|-|:=deg(\omega.-) : N \to \mathbb{N}.$$
We suggest to form the following generating series
$$Z(X,\nu;(t,z)):=exp\left(-\sum_{\beta\in N}\frac{N_\beta}{|\beta|}z^\beta\right),$$
where $N_\beta \in \Ql(t)$ is the total Euler series of $\nu_{1,0}$
on $\RMgn_{1,0}(X,\beta)$. This generating series is an element of $\Ql(t)[[N]]$, the 
completed group ring of $N$ over the field of rational function $\Ql(t)$. 

We like to think of $Z(X,\nu;(t,z))$ as the algebraic analogue of Ruelle's dynamical
zeta function (see e.g.
\cite[\S 1.4-(B)]{Ruelle}) that counts periodic orbit of the flow associated to a global vector field.
It is tempting to ask if $Z(X,\nu;(t,z))$ is rational, but we currently have not enough intuition 
whether this is a reasonable statement to expect or not. We also note that 
similar generating series, without any vector field considerations, already appear
in \cite{MR1659474}. \\

\subsection{Higher genus considerations}\label{sectionhigher}

So far we have started with a global vector field $\nu$ on $X$, so that  
(irreducible) invariant curves can only be of genus $0$ or $1$. In particular, 
the zero locus of the induced vector field $\nu_{g,n}$ on $\RMgn_{g,n}(X,\beta)$ is not a very pertinent 
object to consider when $g>1$, as this corresponds
to stable maps that factor through  the zero locus of $\nu$ in $X$. 
In order to get something interesting
for all genuses, we have to allow ourselves to replace vector fields by more general \emph{foliations}. 

It is a general fact from derived foliation theory, see \cite{tvrh, tvindex}, that 
any foliation $\F$ on $X$ induces a derived foliation $\F_{g,n}$ on 
$\RMgn_{g,n}(X,\beta)$, but again, when $g>1$ the singularities of this
induced foliation is a rather poor geometric object  (it consists essentially of pointwise
invariant curves, i.e. stable maps to the singular locus of the original
foliation on $X$). It is however possible to replace the zero locus of $\nu_{g,n}$ on
$\RMgn_{g,n}(X,\beta)$ with a
more relevant derived stack $\RMgn_{g,n}(X,\beta)^\F$ of 
\emph{$\F$-equivariant stable maps} defined as follows. For simplicity of exposition, we 
suppose that $\F$ is a one dimensional foliation on $X$, defined by 
a line bundle $\mathcal{L}$ together with a morphism of $\OO_X$-modules $s : \mathcal{L} \to T_X$. 
This defines a one dimensional foliation $\F$ on $X$, possibly with singularities at those 
points in $X$ where $s$ vanishes.
By definition, an $S=\Spec\, A$-point in the
derived stack $\RMgn_{g,n}(X,\beta)^\F$ consists of a pair
$(f,u)$, where $f : C \to X$ is stable map (where $C \to S$ is a relative
curve of genus $g$ with n marked points), and $u$ is a morphism of $\OO_C$-modules
$u : f^*(\mathcal{L}) \longrightarrow \T_C(-\Sigma),$
making commutative the diagram 
$$\xymatrix{
f^*(\mathcal{L}) \ar[r]^-{s} \ar[d]_-{u} & f^*(\T_X)  \\
\T_{C/S}(-\Sigma) \ar[r] & T_{C/S}, \ar[u]_-{Tf}
}$$
where $\Sigma=\coprod_n S \hookrightarrow C$ is the relative divisor of marked points in $C$.
This will define
a derived Deligne-Mumford stack $\RMgn_{g,n}(X,\beta)^\F$, endowed with a forgetful map
to $\RMgn_{g,n}(X,\beta)$, from which we can indeed
extract interesting numerical invariants. As in the case studied in this paper, 
these derived stacks are not quasi-smooth, and thus only Euler series as in Section \ref{eulseriessec} will make sense in 
general, 
unless specific conditions on $\F$ ensures that we have quasi-smoothness. When $\mathcal{L}$ is the trivial
line bundle, so $s$ is global vector field on $X$, the derived stack $\RMgn_{g,n}(X,\beta)^\F$
coincides with the derived locus of $\nu_{g,n}$. Essentially all of our results and constructions in the previous sections extend from
vector fields to more general one dimensional foliations, simply by using $\RMgn_{g,n}(X,\beta)^\F$.
The general study of the derived stack of $\F$-equivariant stable maps
will appear elsewhere.

\section{Appendix: K-theory and G-theory of quasi-smooth derived stacks}

We collect in this Appendix a few facts about the $K$-theory and $G$-theory of derived stacks 
that we have used in the main text. The interested reader will find more details e.g. in 
\cite{khank}.\\

Recall that for a $k$-cdga $A$ (non-positively graded and Noetherian, for simplicity), an $A$-dg-module $M$ is said to be \emph{coherent} if it is cohomologically bounded, and each $H^i(M)$ is finitely generated over $H^0(A)$.\\ 
If $X$ is a Noetherian derived Artin stack, we denote by $\QCoh(X)$ the derived $\infty$-category of quasi-coherent complexes on $X$
(see \cite{toenems}, 3.1). $\QCoh(X)$ is a symmetric monoidal $\s$-category (via the derived tensor product over $\mathcal{O}_X)$, and we denote by $\Perf(X)$ its full subcategory consisting of objects which are dualizable with respect to this tensor product. We write $\Coh(X)$ (respectively $\Coh(X)_{\geq 0}$, respectively $\Coh(X)_{\leq 0}$)  for the full subcategory of $\QCoh(X)$ consisting of complexes $E$ such that, for any smooth map $u: S=\mathrm{Spec}\, A \to X$, $u^*E$ is a coherent (respectively, and connective, respectively, and co-connective) $A$-dg module. The pair $(\Coh(X)_{\geq 0},\Coh(X)_{\leq 0})$ defines a t-structure on $\Coh(X)$ (i.e. a usual t-structure on the homotopy category $\mathrm{h}\Coh(X)$, \cite[1.2]{luha}), whose heart will be denoted by $\Coh_{\heartsuit}(X)$. Note that $\Coh_{\heartsuit}(X)$ is equivalent to the abelian category of coherent complexes on the truncation $\tau_0(X)$ of $X$. \\
Both $\Perf(X)$ and $\Coh(X)$ are stable $\s$-categories (\cite[1.1]{luha}), and therefore  (\cite{barw, bgtab}) we can consider their (connective) $K$-theory spectra, $\mathbf{K}(X):= \mathbf{K}(\Perf(X))$, $\mathbf{G}(X):=\mathbf{K}(\Coh(X))$.\\

Recall that a morphism $f: X\to Y$ of derived Artin stacks is (derived) \emph{lci} if it is locally of finite and the relative cotangent complex $\mathbb{L}_{f}$ is of Tor-amplitude $\leq 1$. Then, $X$ is \emph{quasi-smooth} (or, synonymously, \emph{lci}) if $X \to \mathrm{Spec}\, \mathbb{C}$ is lci. Note that an arbitrary (a priori non-noetherian) derived Artin stack that is lci over an arbitrary regular noetherian base derived Artin stack is automatically noetherian and its structure sheaf is coherent.\\ It is easy to see that an lci morphism $X\to Y$ is of finite Tor-amplitude (e.g. \cite[Lemma 1.15]{khank}). The pullback $f^*:\mathbf{K}(Y) \to \mathbf{K}(X)$ exists for any morphism $f: X\to Y$ of derived Artin stacks (because the inverse image preserves perfect complexes), and the pushforward $f_*: \mathbf{K}(X) \to \mathbf{K}(Y)$ exists if $f$ is proper and lci.
On $G$-theory, a morphism $f:X\to Y$ induces a pullback $f^*: \mathbf{G}(Y) \to \mathbf{G}(X)$ if $f$ has finite Tor-dimension (e.g. it is lci) because the inverse image preserves coherence, and $f$ induces a pushforward $f_*: \mathbf{G}(X) \to \mathbf{G}(Y)$ if $f$ is proper and of finite cohomological dimension.\\
The pushforward along the inclusion $\tau_0(X) \to X$ induces an equivalence $\mathbf{G}(\tau_0(X)) \simeq \mathbf{G}(X)$ (\cite[Thm. 6.1]{barw}).\\

Note that, if $\mathcal{E}\in \Perf (X)$ and $\mathcal{F} \in \Coh(X)$, then  $\mathcal{E}\otimes \mathcal{F} \in \Coh (X)$ (derived tensor product), so that we get a well defined map of spectra $$\cap: \mathbf{K}(X)\otimes \mathbf{G}(X) \to \mathbf{G} (X)$$ inducing a $\mathbf{K}(X)$-module structure on $\mathbf{G}(X)$.\\
We have the following \emph{projection formula} (\cite[Proposition 3.7]{khank}). If $f: X\to Y$ is proper and of finite cohomological dimension, we have the following equality in $\mathbf{G}(Y)$ $$y \cap f_*(x)= f_*(f^*(y)\cap x)$$ for any $x \in \mathbf{G}(X)$, and $y \in \mathbf{K}(Y)$.\\

Given a (homotopy) pullback diagram of derived Artin stacks $$\xymatrix{X' \ar[r]^-{g'} \ar[d]_-{f'} & X \ar[d]^-{f} \\
Y' \ar[r]^-{g} & Y}$$ with $g$ proper and of finite cohomological dimension, and $g$ is of finite Tor-amplitude, we have a canonical homotopy $f^*g_* \simeq g'_{*}f'^{*}$ of morphisms $\mathbf{G}(Y') \to \mathbf{G}(X)$ (see \cite[Proposition 3.8]{khank}).\\

For a quasi-smooth derived Artin stack $X$ of finite type (over $k$), there is a Grothendieck-Riemann-Roch transformation $$\tau_{X}: \mathrm{G}_0(X) \longrightarrow \mathrm{H}_{2*}^{BM}(X,\Ql(-*)):=\oplus_{n}\mathrm{H}_{2n}^{BM}(X,\Ql(-n))$$ where $H_{2n}^{BM}(X,\Ql(-n))$ denotes the $2n$-th Borel-Moore homology of the truncation $\tau_0X$ (\cite[6.3]{khank}, \cite[Def. 2.1]{khanvir}).\\

Note that an arbitrary (a priori non-noetherian) derived Artin stack that is lci over an arbitrary regular noetherian base derived Artin stack is automatically noetherian and its structure sheaf is coherent. In particular, a quasi-smooth derived Artin stack $X$ over $k$ has a fundamental class $[X]=[\mathcal{O}_X] \in \mathrm{G}_0(X)$, which under the isomorphism $\mathrm{G}_0(X)\simeq \mathrm{G}_0(\tau_0X)$ corresponds to a class $[X]_{vir}\in \mathrm{G}_0(\tau_0X)$, called the \emph{virtual fundamental class} of $X$. We have $$[X]^{vir}= \sum_{i}(-1)^i[H^i(\mathcal{O}_X)]$$ (the sum being finite since $\mathcal{O}_X$ is bounded). Analogously (\cite[Construction 3.6]{khanvir}), there is a fundamental class $[X]_{\mathrm{BM}} \in \mathrm{H}_{2d}^{BM}(X,\Ql(d)$ in Borel-Moore homology, where $d$ is the virtual dimension of $X$ (i.e. the Euler characteristic of its cotangent complex). The classes $\tau_X([X])$ and $[X]_{\mathrm{BM}}$ differ by multiplication by the inverse of the Todd class of the cotangent complex of $X$ (Kontsevich Formula, see \cite[Thm. 6.12]{khank}).\\
There is a Chern character ring morphism $$Ch : K_0(X) \to H^{2*}(X,\Ql(*)):=\oplus_{n}\mathrm{H}^{2n}(X,\Ql(n))$$ to $\ell$-adic cohomology of $X$, that enters in a Grothendieck-Riemann-Roch (GRR) formula for lci morphisms; we do not need such a GRR formula in the main text. \footnote{The reader may look at \cite[Thm. 6.22]{khank} where a version of the GRR formula is proved (with values in a different cohomology theory).}

\bibliographystyle{alpha}
\bibliography{BrokP1-Bib.bib}

\end{document}